\documentclass{article}[12pt]
\usepackage{amssymb, amsmath, amsthm, verbatim, float}
\usepackage{sectsty}
\usepackage{color, epsfig}

\newcommand\blfootnote[1]{%
  \begingroup
  \renewcommand\thefootnote{}\footnote{#1}%
  \addtocounter{footnote}{-1}%
  \endgroup
  }

\usepackage{amssymb, amsmath, amsthm, verbatim, float}

\usepackage{color, epsfig}

\usepackage{graphicx}
\usepackage{amsmath}
\usepackage{amssymb}
\usepackage{multicol}

\sectionfont{\centering}

\newcommand{\grow}{{\mathcal{G}}}
\newcommand{\coarsen}{{\mathcal{C}}}

\newcommand{\re}{{\mathbb{R}}}
\newcommand{\intgr}{{\mathbb{Z}}}

\newcommand{\cma}{{\hbox{ , }}}
\newcommand{\interp}{{\mathcal{I}}}

\newcommand{\gbold}{{\boldsymbol{g}}}

\newcommand{\lr}[1]{{\langle #1 \rangle}}
\newcommand{\hl}{{h_l}}
\newcommand{\hlp}{{h_{l+1}}}

\newcommand{\om}[1]{{\Omega_{{#1}}}}
\newcommand{\omri}{{\Omega_{R,\ibold}}}
\newcommand{\omria}{{\Omega_{R,\ibold,\alpha}}}
\newcommand{\omrib}{{\Omega_{R,\ibold,\beta}}}

\newcommand{\omriprah}{{\Omega_{R,\ibold',\alpha}^h}}
\newcommand{\omriprbh}{{\Omega_{R,\ibold',\beta}^h}}

\newcommand{\omripraH}{{\Omega_{R,\ibold',\alpha}^H}}
\newcommand{\omriprbH}{{\Omega_{R,\ibold',\beta}^H}}
\newcommand{\omrih}{{\Omega_{R,\ibold}^h}}
\newcommand{\omriah}{{\Omega_{R,\ibold,\alpha}^h}}
\newcommand{\omribh}{{\Omega_{R,\ibold,\beta}^h}}

\newcommand{\omriaH}{{\Omega_{R,\ibold,\alpha}^H}}
\newcommand{\omribH}{{\Omega_{R,\ibold,\beta}^H}}
\newcommand{\omli}{{\om{R_l,\ibold}}}
\newcommand{\omlih}{{\Omega_{R_l,\ibold}^{h_l}}}
\newcommand{\omlia}{{\Omega_{R_l,\ibold,\alpha}}}
\newcommand{\omliab}{{\Omega_{R_l,\ibold,\beta}}}
\newcommand{\omliah}{{\Omega_{R_l,\ibold,\alpha}^{h_l}}}
\newcommand{\omliabh}{{\Omega_{R_l,\ibold,\beta}^{h_l}}}
\newcommand{\omliahp}{{\Omega_{R_{l+1},\ibold',\alpha}^{h_{l+1}}}}
\newcommand{\omliabhp}{{\Omega_{R_{l+1},\ibold',\beta}^{h_{l+1}}}}
\newcommand{\omlihp}{{\Omega_{R_{l+1},\ibold'}^{h_{l+1}}}}
\newcommand{\mlcOne}{MLC-0}\newcommand{\mlcTwo}{MLC}
\newcommand{\loc}{{loc, \xbold}}
\def\L2norm#1{\| #1 \|_2}

\newcommand{\nref}{{{N}_{ref}}}

\newcommand{\xbold}{{\boldsymbol{x}}}

\newcommand{\pbold}{{\boldsymbol{p}}}

\newcommand{\jbold}{{\boldsymbol{j}}}
\newcommand{\sbold}{{\boldsymbol{s}}}

\newcommand{\ibold}{{\boldsymbol{i}}}

\newcommand{\lbold}{{\boldsymbol{l}}}
\newcommand{\Dim}{{\mathbf{D}}}

\newcommand{\beqa}{\begin{eqnarray*}}
\newcommand{\eeqa}{\end{eqnarray*}}

\newcommand{\Domain}[1]{\ifmmode\mbox{\tt Domain<}#1\mbox{\tt>}\else{\tt Domain<$#1$>}\fi}

\newcommand{\fnc}[2]{\ifmmode\mbox{\tt#1(}#2\mbox{\tt)}\else{\tt#1(}$#2${\tt)}\fi}

\newcommand{\IEFab}[1]{\ifmmode\mbox{\tt IEFab<}#1\mbox{\tt>}\else{\tt IEFab<$#1$>}\fi}

\newcommand{\Point}[1]{\ifmmode\mbox{\tt Point<}#1\mbox{\tt>}\else{\tt Point<$#1$>}\fi}

\newcommand{\RectDomain}[1]{\ifmmode\mbox{\tt RectDomain<}#1\mbox{\tt>}\else{\tt RectDomain<$#1$>}\fi}

\newcommand{\sign}[1]{\ifmmode\mbox{sign}(#1)\else sign(#1)\fi}
\newcommand{\Vector}[1]{\ifmmode\mbox{\tt Vector<}#1\mbox{\tt>}\else{\tt Vector<#1>}\fi}

\newcommand{\bi}{\begin{itemize}}
\newcommand{\ei}{\end{itemize}}

\newcommand{\D}{\begin{itemize} \item[]}
\newcommand{\bv}{\begin{verbatim}}
\newcommand{\ev}{\end{verbatim}}

\begin{document}
\title{Computation of Volume Potentials on Structured Grids via the Method of Local Corrections}
\author{Chris Kavouklis, Phillip Colella \\ Computational Research Division \\ Lawrence Berkeley National Laboratory \\ 1 Cyclotron Road, Berkeley, CA 94720, United States}
\maketitle
\begin{abstract}
We present a new version of the Method of Local Corrections (MLC) \cite{mlc}, a multilevel, low communications, non-iterative, domain decomposition algorithm for the numerical solution of the free space Poisson's equation in 3D on locally-structured  grids. In this method, the field is computed as a linear superposition of local fields induced by charges on rectangular patches of size $O(1)$ mesh points, with the global coupling represented by a coarse grid solution using a right-hand side computed from the local solutions. In the present method, the local convolutions are further decomposed into a short-range contribution computed by convolution with the discrete Green's function for an $Q^{th}$-order accurate finite difference approximation to the Laplacian with the full right-hand side on the patch, combined with a longer-range component that is the field induced by the terms up to order $P-1$ of the Legendre expansion of the charge over the patch. This leads to a method with a solution error that has an asymptotic bound of $O(h^P) + O(h^Q) + O(\epsilon h^2) + O(\epsilon)$, where $h$ is the mesh spacing, and $\epsilon$ is the max norm of the charge times a rapidly-decaying function of the radius of the support of the local solutions scaled by $h$. Thus we have eliminated the low-order accuracy of the original method (which corresponds to $P=1$ in the present method) for smooth solutions, while keeping the computational cost per patch nearly the same with that of the original method. Specifically, in addition to the local solves of the original method we only have to compute and communicate the expansion coefficients of local expansions (that is, for instance, 20 scalars per patch for $P=4$). Several numerical examples are presented to illustrate the new method and demonstrate its convergence properties.
\end{abstract}


\section{\large{Introdu]ction}}
We are interested in solving Poisson's equation with infinite domain boundary conditions in three dimensions, that is
\begin{eqnarray}
\label{eqn:poisson}
&&\Delta{\phi} \equiv \frac{\partial^2 \phi}{\partial x^2}+\frac{\partial^2 \phi}{\partial y^2}+\frac{\partial^2 \phi}{\partial z^2}=f,\:\text{in}\:\:\mathbb{R}^3 ,\\
&&\nonumber \phi(\boldsymbol{x})=-\frac{1}{4\pi \|\boldsymbol{x}\|}\int_{\mathbb{R}^3}f(\boldsymbol{y})d\boldsymbol{y}+o\left(\frac{1}{\|\boldsymbol{x}\|}\right),\:\:\|\boldsymbol{x}\|\rightarrow \infty ,
\end{eqnarray}
where $f$ is a function with bounded support and by $\|\cdot\|$ we denote the Euclidean norm. It is well known that problem (\ref{eqn:poisson}) has a solution if $f$ is H{\"o}lder continuous and has compact support $\Omega$ \cite{gilbarg}. Furthermore, the solution of (\ref{eqn:poisson}) is unique by means of a maximum principle argument for harmonic functions and is given as a convolution of the data with the three dimensional infinite domain Green's function \cite{evans}
\begin{gather}
\label{eqn:potential}
\phi (\boldsymbol x)=\int_{\Omega}G(\boldsymbol x - \boldsymbol y)f(\boldsymbol y) d \boldsymbol y \equiv (G*f)(\boldsymbol x), \\
G(\boldsymbol z)=-\frac{1}{4\pi\|\boldsymbol z\|} . \nonumber
\end{gather} 
In addition, if $\:\Omega\subset B(\boldsymbol{x}_0,R)\:$, where $\:B(\boldsymbol{x}_0,R)\:$ is the closed ball of radius $\:R$ centered at point $\boldsymbol{x}_0$, then $\phi$ is harmonic in $\mathbb {R}^{3} \backslash B(\boldsymbol{x}_0,R)$ and hence real analytic. By differentiating (\ref{eqn:potential}), we find that the derivatives of the potential are rapidly-decaying functions of the form
\begin{equation}
\label{eqn:decay}
(\nabla^\pbold \phi)(\boldsymbol{x})=O\left(\left(\frac{1}{||\boldsymbol{x}-\boldsymbol{x}_0||}\right)^{||\pbold||_1+1}R^3||f||_{\infty} \right) .
\end{equation}

This suggests a domain-decomposition strategy, in which the contribution to the fields on each local domain is computed independently and the non-local coupling is computed using a reduced number of computational degrees of freedom.
This approach has been exploited for particle methods with the right hand side in (\ref{eqn:poisson}) given by $f(\boldsymbol{x})=\sum_{i}q_i\delta(\boldsymbol{x}-\boldsymbol{x}_i)$. For instance, we mention the Barnes-Hut algorithm \cite{barneshut}, the Fast-Multipole Method (FMM) \cite{fmmparticle, fmmadaptive2D,fmm3dparticle}, and the Method of Local Corrections (MLC) \cite{anderson,almgrenPhD, almgren}.
The aforementioned particle algorithms have been modified to handle gridded data; for a more comprehensive review that includes benchmark studies of the FFT, FMM and multigrid methods, see \cite{gholami}. 

The present work is based on the extension of the Method of Local Corrections to structured-grid data described in \cite{balls1,balls2,mlc}. In this approach, the support of the right-hand side is discretized with a rectangular grid, which is decomposed into a set of cubic patches. For two levels the method proceeds in three steps: (i) a loop over the fine disjoint patches and the computation of local potentials induced by the charge restricted to those patches on sufficiently large extensions of their support (downward pass); (ii) a global coarse-grid Poisson solve with a right hand side computed by applying the coarse-grid Laplacian to the local potentials of step (i); and (iii) a correction of the local solutions computed in step (i) on the boundaries of the fine disjoint patches based on interpolating the global coarse solution from which the contributions from the local solutions have been subtracted (upward pass). These boundary conditions are propagated into the interior of the patches by performing Dirichlet solves on each patch. This can be generalized by replacing the global coarse solution in (ii) by a recursive call to MLC, or by replacing uniform grids at each level covering the entire domain by nested block-structured locally-refined grids. The local volume potentials are computed using a high-order finite-difference approximation to the Laplacian, combined with an extension to three dimensions of the James-Lackner algorithm \cite{james,lackner} for representing infinite-domain boundary conditions. Furthermore, in order to make the nested refinement version of this algorithm practical, we require that $R=O(H) = O(h)$, where $R$ is the radius (in max norm) of local patches, $H$ the coarse mesh spacing, and $h$ the fine mesh spacing (i.e., a fixed number of points per patch and a fixed refinement ratio). In \cite{mlc}, the local field calculation in (i) was split into two contributions: one that represented the field induced by the complete charge distribution on a patch, and a second corresponding to the monopole component of the charge. By using such a splitting, it is possible to obtain a convergent method by using a relatively large region for computing the monopole component only while keeping the overall computation and communications cost low. However, the convergence properties of the resulting method were erratic, and exhibited a large $O(h)$ solution error for smooth charge distributions that were well-resolved on the fine grid.

In the present work, we generalize the method in \cite{mlc} in a way that preserves the reduced-communication properties of that method, and leads to an error analysis that explains the observed convergence behavior. In particular, we replace the separate treatment of the monopole component of the charge on each patch by a similar treatment of a truncated expansion in Legendre polynomials of the charge distribution on each patch. Our error analysis predicts an $O(h^{P})+O(h^Q)+O(\epsilon h^2) + O(\epsilon)$ solution error, where $P-1$ is the maximum degree of the polynomials in the Legendre expansions, and $Q$ is the order of accuracy of the finite-difference discretization used to compute the local potentials. This is consistent with the earlier results in \cite{mlc} corresponding to $P=1$. The $O(\epsilon)$ term is a localization error, proportional to the max norm of the charge divided by a localization distance (measured multiples of the patch size) raised to the order of accuracy of the discretized Laplacian on harmonic functions. We also change the detailed approach to computing the local potentials, replacing the James-Lackner representation of the infinite--domain boundary conditions in the calculation of the local potentials in step (i) with local discrete convolutions computed using FFTs via a variation on Hockney's domain--doubling method \cite{hockneyMCP}. This leads to a conceptually simpler algorithm, and provides a compact numerical kernel on which to focus the effort of optimization.

In this paper, we focus on the design of the algorithm, including an error analysis of the method and calculations that demonstrate the error properties derived from that analysis. In a second paper \cite{mlc2}, we will present performance and parallel scaling results on high-performance computing platforms.

\section{\large{Mehrstellen Discretization and Finite Difference Localization}}

{\bf Notation}. 
We denote by $D^h,\om{}^h
 \dots \subset \intgr^3$ grids with grid spacing $h$ of discrete points in physical space: $\{\gbold h: \gbold \in D^h\}$. Arrays of values defined over such sets will approximate functions on subsets of $\re^3$, i.e. if $\psi = \psi(\xbold)$ is a function on $D \subset \re^3$, then 
 $\psi^h[\gbold] \approx \psi(\gbold h)$. We denote operators on arrays over grids of mesh spacing $h$ by $L^h, \Delta^h , \dots$;
 $L^h(\phi^h): D^h \rightarrow \re$. Such operators are also defined on functions of $\xbold \in \re^3$, and on arrays defined on finer grids 
 $\phi^{h'}$, $h = N h', N \in \mathbb{N}_+$, by sampling: $L^h(\phi) \equiv L^h(\mathcal{S}^h(\phi))$, $\mathcal{S}^h(\phi)[\gbold] \equiv \phi(\gbold h)$; $L^h(\phi^{h'}) \equiv L^h(\mathcal{S}^h(\phi^{h'}))$, $\mathcal{S}^h(\phi^{h'})[\gbold] \equiv \phi^{h'}[N\gbold]$.
 
 For a rectangle $D=[\boldsymbol{l},\boldsymbol{u}]$, defined by its low and upper corners $\boldsymbol{l},\boldsymbol{u}\in \mathbb{Z}^3$, we define the operators
\begin{gather*}
\mathcal{G}(D,r)=[\boldsymbol{l}-(r,r,r),\boldsymbol{u}+(r,r,r)],r\in\mathbb{Z} \\
\mathcal{C}(D) = \Big[ \Big\lfloor \frac{\boldsymbol{l}}{\nref} \Big \rfloor,\Big\lceil \frac{\boldsymbol{u}}{\nref} \Big \rceil \Big]
\end{gather*} 
Throughout this paper, we will use $N_{ref} = 4$ for the refinement ratio between levels.

 We begin our discussion presenting the finite difference discretizations of (\ref{eqn:poisson}) that we will be using throughout this work and some of their properties that pertain to the Method of Local Corrections. Specifically, we are employing Mehrstellen discretizations \cite{collatz} (also referred to as compact finite difference discretizations) of the 3D Laplace operator
\begin{equation} \label{eqn:discreteLaplacian}
(\Delta^h \phi^h)_{\boldsymbol{g}}=\sum_{\boldsymbol{s}\in[-s,s]^3}a_{\boldsymbol{s}} \phi_{\boldsymbol{g}+\boldsymbol{s}}^h,a_{\boldsymbol{s}}\in \mathbb{R}.
\end{equation}
If $\phi^h$ is defined on $D^h$, then $\Delta^h \phi^h$ is defined on $D^{h,s} \equiv \grow(D^h,-s)$.
 The associated truncation error $\tau^h \equiv (\Delta^h - \Delta)(\phi) = -\Delta^h(\phi^h-\phi)$ for the Mehrstellen discrete Laplace operator is of the form
\begin{gather}
\label{eqn:mehrstellenTau}
\tau^h(\phi) = C_2 h^2 \Delta( \Delta \phi ) + \sum_{q'=2}^{\frac{q}{2}-1}h^{2q'}\mathcal{L}^{2q'}(\Delta \phi)+h^qL^{q+2}(\phi)+O(h^{q+2}),
\end{gather}
 where $q$ is even and $\mathcal{L}^{2q'}$ and $L^{q+2}$ are constant-coefficient differential operators that are homogeneous, i.e. for which all terms are derivatives of order $2q'$ and $q+2$, respectively. For the two operators we will consider here, $C_2=\frac{1}{12}$. In general, the truncation error is $O(h^2)$. However, if $\phi$ is harmonic in a neighborhood of $\boldsymbol{x}$, 
\begin{equation}
\label{eqn:exteriortruncerror}
\tau^h(\phi)(\xbold) = \Delta^h(\phi)(\boldsymbol{x})=h^qL^{q+2}(\phi)(\boldsymbol{x})+O(h^{q+2}) .
\end{equation}
In our numerical test cases we make use of the 19-point ($L_{19}^h$) and 27-point ($L_{27}^h$) Mehrstellen stencils \cite{carey} that are described in the Appendix (Section \ref{sec:L19L27}), for which $q=4$ and $q=6$, respectively. In general, it is possible to define operators for which $s = \lfloor \frac{q}{4} \rfloor $ 
for any even $q$, using higher-order Taylor expansions and repeated applications of the identity 
\begin{gather*}
\frac{\partial^{2 r} \phi}{\partial x_d^{2 r}} = \frac{\partial^{2 r-2}}{\partial x_d^{2 r-2}}(\Delta \phi ) - \sum_{d' \neq d} \frac{\partial^{2 r}}{x_{d'}^{2 r - 2} x_d^2} (\phi) .
\end{gather*}
Since we are primarily concerned with solving the free-space problem, the corresponding discrete problem can be expressed formally as a discrete convolution. 
\begin{gather}
\label{eqn:discreteConvolution}
(G^h*f^h)  = (\Delta^h)^{-1} (f^h) \cma (G^h * f^h)[\gbold] \equiv \sum \limits_{\gbold' \in \mathbb{Z}^3} h^3 G^h[\gbold - \gbold'] f[\gbold']^h
\end{gather}
where the discrete Green's function $G^h[\gbold] = h^{-1} G^{h=1}[\gbold]$ satisfies 

\begin{equation}
\label{eqn:dgf}
(\Delta^{h=1} G^{h=1})[\gbold]  =\left\{ \begin{array}{l} 
  1, \hbox{ if } \gbold = \boldsymbol{0}\\
  0, \hbox{ otherwise} 
\end{array} \right.
\end{equation}
\noindent and
\begin{gather*}
G^{h=1}[\gbold] = -\frac{1}{4 \pi ||\gbold||} +o\Big( \frac{1}{||\gbold || } \Big) \cma \|\gbold\| \rightarrow \infty.
\end{gather*}
We use these conditions to construct approximations to $G^h$ numerically, see the Appendix. For any $n$, we have
\begin{gather*}
\sum \limits_{\gbold \in D} h^3 |G^h[\gbold] |  \leq C \cma C = C(n h) \cma D \subseteq [-n, \dots , n]^3 ,
\end{gather*}
from which it follows that convolution with $G^h$ is max-norm stable on bounded domains, i.e. , 
\begin{gather}
||G^h * f^h||_\infty \leq C' ||f^h||_\infty \cma C' \text{ independent of $f$, $h$ }, \\
\nonumber supp(f^h) \subseteq \Big[ - \Big \lfloor \frac{A}{h}\Big \rfloor, \dots ,\Big \lceil \frac{A}{h} \Big\rceil \Big]^3
\label{eqn:convolutionbound}
\end{gather} 
for any fixed $A > 0$.

The form of the truncation error (\ref{eqn:mehrstellenTau}) allows us to compute $q^{th}$-order accurate solutions to (\ref{eqn:poisson}) by modifying the right-hand side, i.e.
\begin{gather}
\Delta^h (\phi) = \tilde{f}^h + O(h^q) \\
\tilde{f}^h = f^h + \Big(C_2 h^2 (\Delta(f))^h + \sum_{q'=2}^{\frac{q}{2}-1}h^{2q'}\mathcal{L}^{2q'}(f)^h \Big)  ,  \label{eqn:msCorr}
\end{gather}
and replacing the differential operators on the right-hand side with finite difference approximations. If only a fourth-order accurate solution is required, it suffices to use the first term, leading to a correction of a particularly simple form:
\begin{gather}
\label{eqn:msCorr4}
\phi = G^h*f^h + C_2 h^2 f^h + O(h^4).
\end{gather}
In particular, the solution error $\epsilon^h = G^h*f^h - \phi = O(h^4)$ away from the support of $f$ without any modification of $f^h$.

\begin{figure}[t]
\centering
\includegraphics[width=0.49\textwidth]{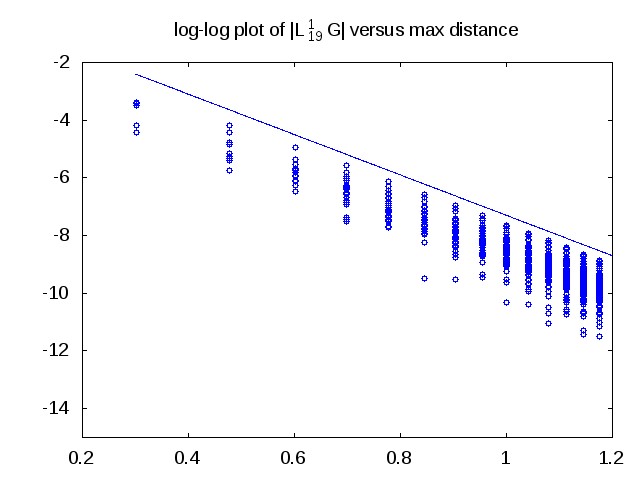}
\includegraphics[width=0.49\textwidth]{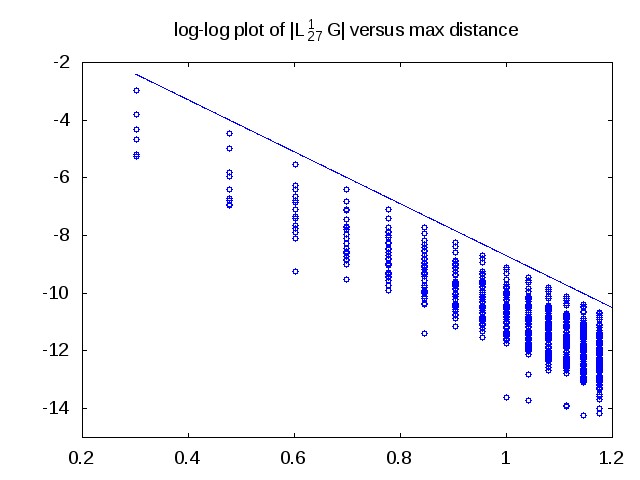}
\caption{Scatter plots of $log_{10}(|(\Delta^{h=1} \mathcal{S}^{h=1}(G))[\gbold]|)$ versus $log_{10}(||\gbold||_\infty)$, $\gbold \in \mathbb{Z}^3$  at points away from the singularity of $G$ for the $L^h_{19}$ and $L^h_{27}$ discrete Laplacians. The slopes of the lines depicted are -7 and -9 for the $L^h_{19}$ and $L^h_{27}$, respectively.}
\label{fig:scatterplots}
\end{figure}

Suppose that $supp(f)\subset P_{\boldsymbol{c}}$, where $P_{\boldsymbol{c}}=\boldsymbol{c}+[-R,R]^3$ is a cube of radius {\it R} centered at point $\boldsymbol{c}$ and that the differential operator $L^q$ is a linear combination of derivatives of order {\it q}. By differentiating (\ref{eqn:potential}), we have
\begin{equation}
\label{eqn:l2} \left[(L^q G)*f\right](\boldsymbol{x})=O\left(\left(\frac{1}{R}\right)^{q-2} \frac{1}{\left\| \frac{\boldsymbol{x}}{R} - \frac{\boldsymbol{c}}{R} \right\|^{q+1}_{\infty}}\right)\|f\|_{\infty}.
\end{equation}
In particular, 
away from the support of $f$, \eqref{eqn:mehrstellenTau} becomes
\begin{align}
\label{eqn:q1} \tau^h(f)= \Delta^h(G*f)(\boldsymbol{x})=  O \left( \left( \frac{h}{R} \right)^q \frac{1}{\left\| \frac{\boldsymbol{x}}{R} - \frac{\boldsymbol{c}}{R}\right\|_\infty^{q+3}} \right)\|f\|_{\infty} .
\end{align}
It is precisely this rapid decay of the truncation error, a consequence of the fact that the local potentials are harmonic away from the supports of the associated charges, that allows us to use a coarse mesh for the global coupling computation. In Figure \ref{fig:scatterplots}, scatter plots of the truncation error for the case of a point charge located at the origin using the 19-point and 27-point Laplacians are depicted. The rapid decay of the truncation error in the far-field and the faster decay with increasing $q$ are evident.
Using this localization property of the Mehrstellen operators, we can reduce the cost of computing the potential (\ref{eqn:potential})  induced by a localized charged distribution to the cost of computing the potential near the support of the charge, using the finite difference localization approach originally introduced in \cite{mayo}. 
We assume that the support of $f$ is contained in cube $D$ of radius $R$ centered at $\boldsymbol{c}$.  First, we compute $\phi = G*f$ in the extended cube $D_{\beta}$ of radius $\beta R,\beta>1$. Then we compute $\phi^H = G^H*F^H$ on $\Omega^H$. The coarse right hand side is defined by:
\begin{equation}
\label{eqn:FDL}
F^H=\left\{\begin{array}{l}
\Delta^H(\phi) \text{ , on } D_\beta^{H,s} \cma D_\beta^{H,s} = \grow(\coarsen (D_{\beta}^H),-s) \\
0 \text{ , on } \Omega^H\setminus D_{\beta}^{H,s}.
\end{array} \right.
\end{equation} 
Using (\ref{eqn:q1}), we have
\begin{align}
\Delta^H(\phi^H-G*f) = & \text{ } 0 \text{ on } D_\beta^{H,s} \nonumber \\
= & \text{ } O\Big(\Big(\frac{H}{R}\Big)^q \frac{1}{(k + \beta)^{q+3}} || f ||_\infty \Big) 
\text{ on } \{\gbold: ((k + \beta) + 1)R \ge \| \gbold H \|_\infty \ge (k + \beta) R \}. \label{eqn:loc1}
\end{align}
\noindent where $\:k \in \mathbb{N}$. One can decompose the annular region $ \{\gbold: ((k + \beta) + 1)R \ge \| \gbold H \|_\infty \ge (k + \beta) R \}$ into $O((k + \beta)^2)$ rectangles, each of which of radius $\leq R$, leading to an analogous decomposition of the right-hand side of \eqref{eqn:loc1} into a sum of terms, each of which is supported on one such rectangle. Applying convolution with $G^H$ to both sides of \eqref{eqn:loc1} represented in terms of such sums leads to a solution error given by
\begin{align}
\phi^H-G*f  =& \sum_{k =0}^\infty \label{eqn:loc2}
 O\Big(\Big(\frac{H}{R}\Big)^q \frac{1}{{(k + \beta)}^{q+3} } || f ||_\infty \Big)\\
 = &  O\Big(\Big(\frac{H}{R}\Big)^q \frac{1}{\beta^{q}} || f ||_\infty \Big) \label{eqn:loc2b}
\end{align}
Thus the accuracy of the potential away from the support of the charge can be improved by decreasing the ratio $H/R$; or, for fixed values of that ratio, by adjusting
$\beta$ or $q$. In any case, the error is only weakly dependent on $f$. In this context, we will refer to $\beta$ as a {\it localization radius}. In addition, \eqref{eqn:loc2b} is true independent of whether or not the right-hand side is modified using the Mehrstellen correction (\ref{eqn:msCorr}). The MLC algorithm combines finite difference localization with domain decomposition into a collection of rectangular patches of size $R$ to obtain a low-communication method for computing volume potentials. In that case, we want to keep the number of mesh points per patch fixed, which leads to \eqref{eqn:loc2} being an $O(1)$ error relative to the mesh spacing. Ultimately, that error is controlled by increasing $\beta$, combined with choosing a discretization with a larger $q$. However, the cost of computing the local convolution $G*f$ on $D^{H,s}_\beta$ scales like $\beta^3$. To reduce that cost, we introduce a second localization radius $\alpha$, $\alpha < \beta$. On $D^{H,s}_\alpha$, we use the full convolution to compute $F^H$. In the remaining annular region, we use a reduced representation based on the field induced by the first few moments of the Legendre expansion of $f$, which is much less expensive to compute.

\section{{Method of Local Corrections  -  Semi-Discrete Case}}

\begin{figure}
\begin{center}
\scalebox{.5}{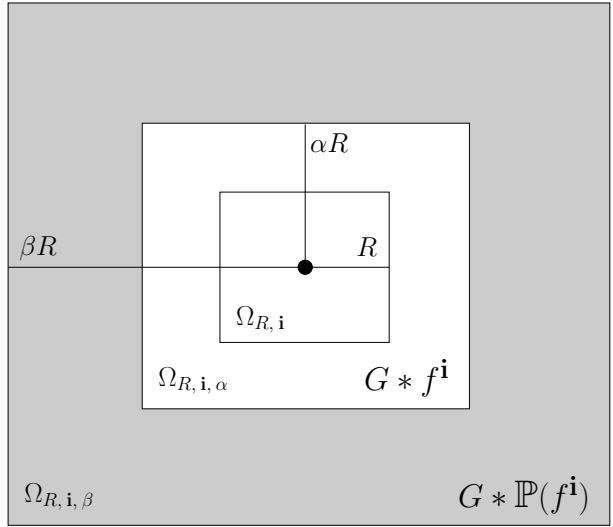}
\caption{Regions associated with subdomain $\omri$. The potential in $\omria$ (white region) is given by $G*f^{\ibold}$. In the ring $\omrib \setminus \omria$ (shaded region) we use the field induced by a truncated Legendre expansion on $\omri$ of the local charge $f^\ibold$ to represent the potential.}
\label{fig:regions1}
\end{center}
\end{figure}

To clarify ideas, we discuss in this section a theoretical proxy for the fully discrete algorithm. We construct a function $\phi^{MLC}:\Omega \rightarrow \re$ that approximates the potential $\phi$ by a linear superposition of local potentials, combined with data interpolated from a discrete global solution. The computational domain is a cube $\Omega$ that contains the support of $f$ and is decomposed into a finite union of disjoint cubic subdomains of equal volume that are translations of $[-R,\:R]^3,R>0$.
\begin{equation}
supp(f)\subset \Omega = \bigcup_{\ibold} \omri \text{ , } \omri=c^{\ibold}+[-R,R]^3 \cma\ibold\in \mathbb{Z}^3 \cma c^\ibold = (2 \ibold + (1,1,1)) R.
\end{equation}
Then $f=\sum_{\ibold} f^\ibold$ where $f^\ibold=f \chi^\ibold$, where $\chi^\ibold$ is the characteristic function of $\omri$. As a consequence, the global potential may be written as
\begin{equation}
\label{eqn:ctsDD}
\phi (\boldsymbol{x})=(G*f)(\boldsymbol{x})=\sum_{\ibold} (G*f^{\ibold})(\boldsymbol{x}).
\end{equation}
In other words, it is the linear superposition of the potentials induced by the local charges $f^\ibold$ which can be computed independently in parallel. The MLC algorithm replaces each of the summands in  (\ref{eqn:ctsDD}) with a solution truncated to zero outside of a  localization radius $\beta R$, with the contribution to the solution outside the localization radius represented by interpolation from a single coarse grid solution $\phi^H$ obtained by summing contributions of the form (\ref{eqn:FDL}) over all the patches. At each point in space, the coarse grid values used to interpolate the global contibution are corrected by subtracting off the contributions of the patches within the localization radius.  
Finally, to reduce the cost of computing the localized potentials, while keeping $\beta$ large enough to make the $O(1)$ contribution to the error coming from localization be acceptably small, we introduce an inner radius $\alpha < \beta$ (see Figure \ref{fig:regions1}). Within that inner radius, we compute the full convolution $G*f^\ibold$; in the annular region  $\omrib \setminus \omria$, the local solution is approximated by $G*\mathbb{P} (f^\ibold )$, where  $\mathbb{P} (f^\ibold)$ is the orthogonal projection onto the Legendre polynomials on $\omri$ of some degree $P-1$:
\begin{gather}
\mathbb{P}(f^\ibold) = \sum \limits_{\pbold \in \mathbb{N}^3: ||\pbold ||_1 < P} \lr{Q^\pbold,f^\ibold} Q^\pbold ,\\
\nonumber Q^{\boldsymbol{p}}(\boldsymbol{x}) 
\label{eqn:led3d}  = R^{-\frac{3}{2}} \prod \limits_{d=1}^3 Q^{p_d}\Big(\frac{x_d - c_d^\ibold}{R}\Big) \cma \boldsymbol{x} \in \omri \cma \boldsymbol{q} \in \mathbb{N}^3,
\end{gather}
where $\lr{\cdot,\cdot}$ is the inner product on $\omri$, and $Q^p:[-1,1] \rightarrow \re$ is the classical Legendre polynomial of degree $p$.

\subsection{The Semi-Discrete MLC Algorithm}
\label{sec:semiDiscreteAlgorithm}

The semi-discrete MLC algorithm consists of three steps. \\

\noindent {\bf Step 1 - Local Convolutions}.\\
We perform local convolutions in regions around each subdomain $\omri$ that are used to compute local charges at points on the grid.

\begin{equation*}
F^{\ibold,H}[\boldsymbol{g}] =\left\{\begin{array}{l}
 \Delta^H(G*f^{\ibold})[\boldsymbol{g}] \text{ , } \boldsymbol{g} \in \omriaH \\\\
 \Delta^H (G*\mathbb{P}(f^{\ibold}))[\boldsymbol{g}]\text{ , } \boldsymbol{g} \in \omribH \setminus  \omriaH  \\\\
 0 \text{ , otherwise}
\end{array} \right.
 \end{equation*}

 \noindent {\bf Step 2 - Global Coarse Solve}. \\
The global charge at coarse mesh points is constructed by assembling local contributions
\begin{equation*}
F^H[\boldsymbol{g}]=\sum_{\ibold}F^{\ibold,H}[\boldsymbol{g}],
\end{equation*}
and we obtain a global approximation $\phi^{H}$ of the potential, represented on the coarse mesh, by computing the discrete convolution over $\Omega^H$.
\begin{equation}
\label{eqn:global_coarse_discr}
\phi^{H} =G^H*F^H.
\end{equation}

\noindent {\bf Step 3 - Local Interactions / Local Corrections.}\\
In the final step, we represent the solution on the boundary of each $\omri$ as the sum of local convolutions induced by charges on nearby patches and values interpolated from the grid calculation, from which the local convolution values have been subtracted. 
\begin{equation}
\label{eqn:mlcdiscr}
\phi^{B,\ibold}(\boldsymbol{x})=\phi^\loc(\boldsymbol{x})+\interp^H(\phi^{H}-\phi^\loc)(\boldsymbol{x}) 
\end{equation}
Here $\interp^H(\psi^H)(\xbold)$  is an interpolation operator that takes as input values of  $\psi^H:\mathcal{N}(\xbold) \rightarrow \re$, where $\mathcal{N}(\xbold) \subset \{\gbold H : \gbold \in \mathbb{Z}^3\}$ and returns a $q_I^{th}$-order accurate polynomial interpolant. In all of the algorithms described here, $\xbold$ and all of the points in $\mathcal{N}(\xbold)$ are coplanar, so the interpolant is particularly easy to construct. 
$\phi^\loc(\xbold)$ is the sum of all local convolutions the support of whose charges is sufficiently close to $\xbold$ so that they contributed to the right-hand side for the grid solution near that point.
\begin{equation}
\label{eqn:philocmlcdiscr}
\phi^\loc(\boldsymbol{x}')=
\sum_{\ibold:\boldsymbol{x} \in \omria}(G*f^{\ibold})(\boldsymbol{x}') +\sum_{\ibold:\boldsymbol{x}\in \omrib\setminus\omria} (G*\mathbb{P}(f^{\ibold}))(\boldsymbol{x}') .
\end{equation}
Equation (\ref{eqn:mlcdiscr}) can be interpreted as the decomposition of the potential at a point $\boldsymbol{x}$, into the sum of local contributions to the potential given by $\phi^\loc$ and corrections to include the global coupling by interpolating a corrected form of the coarse mesh global solution $\phi^{H}$.  Specifically, the correction term in (\ref{eqn:mlcdiscr}) is computed by evaluating $\phi^\loc$ at the points of the interpolation stencil $\mathcal{N}(\boldsymbol{x})$, subtracting these values from $\phi^{H}$ and interpolating the result to $\boldsymbol{x}$.  The MLC solution $\phi^{MLC}$ is specified in terms of solutions to Dirichlet problems on each $\omri$.
\begin{gather}
\label{eqn:last_step}
 \Delta \phi^{MLC}=f^{\ibold} \text{ in }\omri ,  \\
\phi^{MLC}=\phi^{B,\ibold} \text{ on }\partial \omri . \nonumber
\end{gather}

\subsection{Error Analysis} \label{sec:errSemi}
The error of the local corrections step for $\boldsymbol{x}\in\partial \omri$ is given by:
\begin{eqnarray}
\label{eqn:localcorrectionserror}
\nonumber(\phi^{B,\ibold}-\phi)(\boldsymbol{x})&=&\phi^\loc(\boldsymbol{x})-\phi(\boldsymbol{x})-\interp^H\left(\phi^\loc-\phi \right)(\boldsymbol{x})+\interp^H\left(\phi^{H}-\phi\right)(\boldsymbol{x})\\
&=&\epsilon_I^H(\phi^\loc-\phi)(\boldsymbol{x})+\interp^H\left(\phi^{H}-\phi\right) (\boldsymbol{x})
\end{eqnarray}
where $\epsilon_I^H(\psi)(\xbold)$ is the error in applying the interpolation operator $\interp^H$ to a smooth function $\psi$ evaluated on the grid and evaluating it at $\xbold$. There are two sources of error for the semi-discrete algorithm: one from the calculation of $\phi^{H}$ in (\ref{eqn:global_coarse_discr}), and the other due to interpolation at the local corrections step (\ref{eqn:mlcdiscr}). To estimate the former, i.e. the second term of (\ref{eqn:localcorrectionserror}), it suffices to bound the coarse mesh error $\phi^{H}-\phi$. To do so, we estimate the truncation error of the coarse solve (\ref{eqn:global_coarse_discr}) at points $\boldsymbol{g}$:
\begin{eqnarray}
\nonumber \tau_C^H&=&\Delta^H(\phi^{H}-\phi)[\boldsymbol{g}]\\
\label{eqn:SDC}&=&-\Delta^H\left(\sum_{\ibold: \gbold H \notin \omrib}G*f^\ibold \right)[\boldsymbol{g}]
- \Delta^H\left(\sum_{{\ibold:\gbold H \in \omrib \setminus \omria} }G*\left((\mathbb{I}-\mathbb{P})(f^{\ibold})\right)\right)[\boldsymbol{g}]
\end{eqnarray} 
To bound the first term of (\ref{eqn:SDC}), we use (\ref{eqn:q1}) to find that
\begin{align}
\nonumber \Delta^H\left( \sum_{\ibold: \gbold H \notin \omrib} G*f^{\ibold} \right)[\boldsymbol{g}] =&O\left( \left(\frac{H}{R} \right)^q  \sum \limits_{k = 0}^\infty ~~\sum \limits_{\ibold: \gbold H \in \Omega_{R,\ibold,\beta+k+1}\setminus\Omega_{R,\ibold,\beta+k}} \frac{1}{(\beta + k)^{q+3}}  \| f^\ibold \|_{\infty} \right)\\
 =&O\left( \left(\frac{H}{R}\right)^q  \frac{1}{\beta^q}  \|f \|_{\infty} \right) . \label{eqn:estim2discr}
\end{align}
The second term of (\ref{eqn:SDC}) is bounded in a similar fashion.
\begin{align}
\label{eqn:estim2discr2} 
\Delta^H\left( \sum_{{\ibold:\gbold H \in \omrib \setminus\omria} }G*\left((\mathbb{I}-\mathbb{P})(f^{\ibold})\right)    \right)[\boldsymbol{g}] = & O\Bigg( \left(\frac{H}{R}\right)^q \frac{1}{\alpha^q} max_{\ibold} || (\mathbb{I}-\mathbb{P})(f^{\ibold}) ||_\infty \Bigg)
\\ 
= & O\left( \left(\frac{H}{R}\right)^{q} \frac{1}{\alpha^q} H^{P} \right), \nonumber
\end{align}
where we have used
\begin{equation}
\\||(\mathbb{I}-\mathbb{P})(f^{\ibold})\|_{\infty}=O\left(R^{P}\right),
\end{equation}
which follows directly from Taylor's theorem for $f^{\ibold}$ and the fact that $\pi=\mathbb{P}(\pi)$ for polynomials $\pi$ of degree less than $P$. 
As a result, the following estimate for the coarse mesh error holds
\begin{equation}
\label{eqn:coarsemesherrorestimate}
\Delta^H(\phi^{H}-\phi)=O \left( \left(\frac{H}{R}\right)^{q} \frac{1}{\alpha^q} H^{P} \right)+O\left( \left(\frac{H}{R}\right)^q  \frac{1}{\beta^q} \|f\|_{\infty} \right)
\end{equation}
uniformly on coarse mesh points. Since convolution with $G^H$ and the interpolation operator $\interp^H$ are max-norm bounded, $\epsilon_C^H \equiv \phi^{H}-\phi$ is also bounded by an expression of the form of the right-hand side of  (\ref{eqn:coarsemesherrorestimate}).

To bound the first term in (\ref{eqn:localcorrectionserror}), it follows from the fact that the interpolation method is $q_I^{th}$-order accurate that
\begin{align}
\nonumber \epsilon_I^H(\phi^\loc-\phi) (\xbold )=&  H^{q_I}L^{q_I}_I(\phi^\loc-\phi)(\boldsymbol{\xi}) \\
\label{eqn:epsilonIH}=& -H^{q_I} \left( \sum_{{\ibold:\xbold \in \omrib \setminus \omria} }((L^{q_I}_IG)*(\mathbb{I}-\mathbb{P})(f^\ibold))(\boldsymbol{\xi}) + \sum_{\ibold: \xbold \notin \omrib}((L^{q_I}_IG)*f^\ibold)(\boldsymbol{\xi}) \right)
\end{align}
where $\boldsymbol{\xi}$ is in an $O(H)$ neighborhood of $\mathcal{N}(\boldsymbol{x})$ and $L^{q_I}_I$ is a linear differential operator with terms that are derivatives of order $q_I$. Using (\ref{eqn:l2}), a similar argument to that given in the proof of (\ref{eqn:coarsemesherrorestimate}) leads to:
\begin{equation*}
\epsilon_I^H=H^{P+2}O\left(\left(\frac{H}{R}\right)^{q_I-2}\frac{1}{\alpha^{q_I-2}} \right)+H^2O\left(\left(\frac{H}{R}\right)^{q_I-2}\frac{1}{\beta^{q_I-2}}\|f\|_{\infty} \right)
\end{equation*}
so that (\ref{eqn:localcorrectionserror}) is estimated as
\begin{eqnarray}
\label{}
\nonumber \epsilon^{SD} \equiv \phi^{B,\ibold}-\phi&=&H^{P+2}O\left(\left(\frac{H}{R}\right)^{q_I-2} \frac{1}{\alpha^{q_I-2}} \right)+H^2O\left(\left(\frac{H}{R}\right)^{q_I-2}\frac{1}{\beta^{q_I-2}}\|f\|_{\infty} \right)\\
\label{eqn:finemesherrorMLCIIdiscr}&+&O \left( \left(\frac{H}{R}\right)^{q} \frac{1}{\alpha^q} H^P \right)+O\left( \left(\frac{H}{R}\right)^q  \frac{\|f\|_{\infty}}{\beta^{q}} \right)
\end{eqnarray}

\section{\large{Method of Local Corrections  -  Fully-Discrete Case}}
\label{sec:fullyDiscreteAlgorithm}

In this section, we describe the two-level algorithm as it is actually implemented. $\Omega^h$ is a fine-grid discretization of a bounded domain $\Omega$, the latter containing the support of $f$. 
$\Omega^h$ is assumed to be a finite union of rectangles of the form $\omrih = n \ibold + [0,n]^3$, $R = n h / 2$. We also define discrete forms of $\omriah$, $\omribh$:  $\omriah=\mathcal{G}(\omrih, \lceil\frac{(\alpha-1)n}{2}\rceil )$ and $\omribh=\mathcal{G}(\omrih, \lceil \frac{(\beta-1)n}{2}\rceil )$. 
The coarse grid $\Omega^H$ is assumed to cover all of the fine patch data required for the algorithm described below:
$\grow(\coarsen(\omribh), b) \subset \Omega^H$ where $b$ is the radius of the stencil for the interpolation function $\interp^H$.
We also define a discretized form of the characteristic function of a rectangular patch $D \subset \intgr^3$
\begin{equation*}
\chi_{D}(\boldsymbol{x})=\left\{ 
\begin{array}{l}
\frac{1}{8},\text{ if $\boldsymbol{g}$ is a corner of {\it D}} \\\\
\frac{1}{4},\text{ if $\boldsymbol{g}$ lies on an edge of {\it D}} \\\\
\frac{1}{2},\text{ if $\boldsymbol{g}$ lies on a face of {\it D}} \\\\
1,\text{  if $\boldsymbol{g}$ lies in the interior of {\it D}}\\\\
0,\text{  elsewhere} 
\end{array} \right.
\end{equation*} 
In the fully-discrete algorithm, we replace the local convolutions with local discrete convolutions, e.g. 
$G*f^\ibold \rightarrow G^h*f^{\ibold,h}$, $f^{\ibold,h} = \chi_{\Omega^h_{R,\ibold}} f$, and we take $H = N_{ref}h $. 
\subsection{The Fully-Discrete Two-Level Algorithm}

\noindent {\bf 1.  Step 1 - Local Convolutions.}\\
For each $\omrih$, we compute the potential induced by $f^{\ibold,h} = \chi_{\Omega_{R,\ibold}^h} f^h  $.
\begin{equation}
\label{eqn:initiallocalproblemsmlcII}
\phi^{\ibold,h}=G^h * f^{\ibold,h} \text{ on }\mathcal{G}(\omriah ,\nref b).
\end{equation} 
The Legendre expansion coefficients of $f^{\ibold,h}$ required to compute $\mathbb{P}(f^\ibold)$ are computed with composite numerical integration. We employ Boole's rule if $f$ is given only at points of $\Omega^h$ or Gauss integration if $f$ is specified analytically.
For each $\omrih$ we also compute the associated local charges

\begin{equation} \label{eqn:FDLD}
F^{\ibold,H}[\boldsymbol{g}]=\left\{\begin{array}{l} 
\Delta^H \phi_{\ibold}^h[\boldsymbol{g}] \cma \boldsymbol{g} \in \mathcal{C}(\omriah) \\\\
\Delta^H (G^h*\mathbb{P}^h(f^{\ibold,h}))[\boldsymbol{g}] \cma \boldsymbol{g} \in \mathcal{C}(\omribh) \setminus \coarsen(\omriah) \\\\
0  \cma \boldsymbol{g} \notin  \coarsen(\omribh)
\end{array} \right.
\end{equation}

\noindent The values of $\Delta^H (G^h*Q^{\pbold})$ can be computed once and stored, reducing the calculation of $\Delta^H(G^h*\mathbb{P}^h (f^{\ibold,h}))$ to computing linear combinations of the appropriate subset of those precomputed values.\\

\noindent {\bf 2. Global Coarse Solve.}\\
\begin{gather*}
\phi^H = G^H*F^H \hbox{ on } \Omega^H  \cma F^H=\sum_{\ibold} F^{\ibold,H} .
\end{gather*}

\noindent {\bf 3. Local Interactions - Local Corrections}.\\
We define the local potentials at fine boundary points $\boldsymbol{g} \in \partial \omrih$ as combinations of short-range and intermediate-range components
\begin{equation}
\label{eqn:philocmlcII}
\phi^{loc,\gbold}[\boldsymbol{g'}]=\sum_{\ibold':\boldsymbol{g}\in\omriprah}\phi^{\ibold',h}[\boldsymbol{g'}]+\sum_{\ibold':\boldsymbol{g} \in\omriprbh \setminus\omriprah} (G^h*\mathbb{P}^h(f^{\ibold',h}))[\boldsymbol{g'}] , 
\end{equation}
and we correct them by adding the far-field effects as in (\ref{eqn:mlcdiscr})
\begin{equation}
\label{eqn:mlcfulldiscr}
\phi^{B,\ibold, h}[\boldsymbol{g}]=\phi^{loc,\gbold}[\boldsymbol{g}]+\interp^H\left(\phi^{H}-(\phi^{loc,\gbold}) \right)(\boldsymbol{g}h) , \gbold \in \partial \Omega_\ibold^h.
\end{equation}
The interpolation operator on coplanar points $\interp^H$ that we are employing is the same as in \cite{mlc}. 
Using these boundary conditions, we solve the following local Dirichlet problems on $\Omega_{\ibold}^h$ patches
\begin{align}
\label{eqn:last_step_fulldiscr}
\Delta^h \tilde{\phi}^{MLC,\ibold,h}=&f^{\ibold,h} \text{ on } \omrih - \partial \omrih , \\
 \tilde{\phi}^{MLC,\ibold,h}=&\phi^{B,\ibold,h} \text{ on }\partial \omrih.
\nonumber
\end{align}
Finally, the fourth-order Mehrstellen correction (\ref{eqn:msCorr4}) is applied to obtain the values of $\phi^{MLC,h}$
\begin{equation}
\label{eqn:mcapplicationfulldiscr}
\phi^{MLC,h}[\boldsymbol{g}]=\tilde{\phi}^{MLC,\ibold,h}[\boldsymbol{g}]+C_2 h^2 f^h[\boldsymbol{g}] \cma \boldsymbol{g} \in \omrih .
\end{equation}
If we want to go to higher than fourth order accuracy in $h$, the algorithm is more complicated -- the Mehrstellen correction must be applied earlier in the process. We will not discuss the details in this paper.
\subsection{Error Analysis}
We proceed in this section with estimating the error for the fully-discrete MLC algorithm. We want to get some idea of the impact of replacing the analytic continuous convolutions by the discretized convolutions. To do this, we use a modified equation approach, in which we assume that we can approximate the solution error by the action of the operator on the truncation error. In the present setting, this amounts to making the substitution
\begin{gather}
\label{eqn:convMehrst}
G^h * \psi^h  ~ \rightarrow ~ G * (\psi + \delta \tau^h(\psi)) - C_2 h^2 \psi \\
\delta \tau^h(\psi) = \Delta (G^h* \psi^h) - \psi + C_2 h^2 \Delta \psi = O(h^4)
\end{gather}
As in the semi-discrete case, we want to estimate the error in the boundary conditions
\begin{align}
\phi^{B,\ibold,h}[\gbold] - \tilde{\phi}(\gbold h)  = & \phi^{loc,\gbold}[\gbold] - \tilde{\phi}(\gbold h ) + \interp^H(\phi^{H} - \phi^{loc,\gbold})( \gbold h) \nonumber
\\ = & \interp^H(\phi^{H} - \tilde{\phi})(\gbold h) \label{eqn:discreteCoarseError}
\\ + & \phi^{loc,\gbold}[\gbold] - \tilde{\phi}(\gbold h) - \interp^H(\phi^{loc,\gbold} - \tilde{\phi})(\gbold h),\gbold \in \omrih  \label{eqn:discreteInterpError}
\end{align}
where 
\begin{gather*}
\tilde{\phi} \equiv \phi + C_2{h^2}f . 
\end{gather*}

An estimate of the contribution from (\ref{eqn:discreteCoarseError}) is obtained by bounding $\Delta^H (\phi^{H} - \tilde{\phi})$, since $\interp^H$ and convolution with $G^H$ are both stable in max norm. We have, by  \eqref{eqn:convMehrst}, 
\begin{align}
\Delta^H (\phi^{H} - \tilde{\phi})[\gbold] = -  &  \sum \limits_{\ibold':\boldsymbol{g}\notin\omriprbH} \Delta^H(G*f^{\ibold'})[\gbold] - \sum \limits_{\ibold':\boldsymbol{g} \in\omriprbH \setminus\omripraH} \Delta^H(G*(\mathbb{I} - \mathbb{P})(f^{\ibold'}) )[\gbold] \label{eqn:FDterm1}\\ 
- & \sum \limits_{\ibold':\boldsymbol{g}\notin\omriprbH} \Delta^H (G*(\delta \tau^h(f^{\ibold',h})))[\gbold] \nonumber \\
 - & \sum \limits_{\ibold':\boldsymbol{g} \in\omriprbH\setminus\omripraH} \Delta^H (G*\delta\tau^h((\mathbb{I} - \mathbb{P}) (f^{\ibold'})))[\gbold]\nonumber\\
 - & \sum \limits_{\ibold':\boldsymbol{g} \in\omriprbH \setminus\omripraH} \Delta^H (G^h*((\mathbb{P} (f^{\ibold'}))^h - (\mathbb{P}^h (f^{\ibold',h}))))[\gbold] \nonumber + O(h^4)
\end{align}
The first two terms are identical to the ones that appear in the semi-discrete case, while \eqref{eqn:convMehrst}, and the estimate $||(\mathbb{P} - \mathbb{P}^h)(f^{\ibold'})||_\infty=O(h^6)$ (which holds since our quadrature rules for computing the Legendre coefficients are at least sixth-order accurate) guarantee that the remaining terms are $O(h^4)$ or smaller.
Using similar arguments to those in  (\ref{eqn:FDterm1}), we have
\begin{align*}
\phi^{loc,\gbold} - \tilde{\phi} =  -\sum \limits_{\ibold':\boldsymbol{g}\notin\omriprbh} G*f^{\ibold'} - \sum \limits_{\ibold':\boldsymbol{g} \in\omriprbh\setminus\omriprah } G*((\mathbb{I} - \mathbb{P})(f^{\ibold'})) + O(h^4),
\end{align*}
%
and therefore, following \eqref{eqn:epsilonIH}, we have
\begin{align*}
\epsilon_I^H(\phi^{loc,\gbold} - \tilde{\phi})(\gbold h)  =  & H^{P+2}O\left(\left(\frac{H}{R}\right)^{q_I-2}\frac{1}{\alpha^{q_I-2}} \right)+H^2O\left(\left(\frac{H}{R}\right)^{q_I-2}\frac{1}{\beta^{q_I-2}}\|f\|_{\infty} \right) + O(h^4) ,
\end{align*}

Thus we have 
\begin{gather*}
\phi^{B,\ibold,h}[\gbold] - \tilde{\phi}(\gbold h) = \epsilon^{SD} + O(h^4).
\end{gather*}
The stability of the discretized boundary-value problem implies $\|\phi^{MLC,h}-\phi\|_{\infty} = O( \|\phi^{B,h} - \phi\|_{\infty}) + O(h^4)$, so we finally have the following estimate
\begin{eqnarray}
\nonumber\phi^{MLC,h} - \phi&=&O(h^4)+\epsilon^{SD}\\
\nonumber&=&O(h^4)+H^{P+2}O\left(\left(\frac{H}{R}\right)^{q_I-2}\frac{1}{\alpha^{q_I-2}} \right)+H^2O\left(\left(\frac{H}{R}\right)^{q_I-2}\frac{\|f\|_{\infty}}{\beta^{q_I-2}} \right)\\
\label{eqn:finemesherrorMLC}&+&O \left( \left(\frac{H}{R}\right)^q \frac{1}{\alpha^q} H^P \right)+O\left( \left(\frac{H}{R}\right)^q  \frac{\|f\|_{\infty}}{\beta^{q}} \right) .
\end{eqnarray}
at all fine grid points. 
This error can be written in the form
\begin{gather} 
\phi^{MLC,h} = \phi + O(h^4) + O(h^P)  + O\Big( h^2 ||f||_\infty \frac{1}{\beta^{q_I-2}} \Big) + O\Big( ||f||_\infty  \frac{1}{\beta^{q}} \Big).
\end{gather}
Thus MLC differs from classical finite-difference methods in that there is a contribution to the error that does not vanish as $h \rightarrow 0$, i.e. the right-most summand in \eqref{eqn:finemesherrorMLC}. We refer to this contribution to the error as the {\it barrier error}. 
Note that if we take $q_I = q+2$, we obtain the form of the error given in the Introduction. We have specialized this algorithm to the case of fourth-order accuracy, primarily because it allows us the simplification of applying the Mehrstellen correction (\ref{eqn:mcapplicationfulldiscr}) at the end of the calculation. However, this analysis suggests that, even with this simplification, there might be an advantage to using discretizations of the Laplacian with larger $q$, i.e. ones that are higher order accurate when applied to harmonic functions, since the barrier error is proportional to 
$\beta^{-q}$. We observe this to be the case in the results in Section \ref{sec:examples}.

\section{\large{Multilevel Method of Local Corrections}}
\label{sec:multilevel}
Following \cite{mlc}, we generalize the method in Section \ref{sec:fullyDiscreteAlgorithm} to the case of an arbitrary number of levels $l=0,\dots,l_{max}$, where $l_{max}$ is the finest level on which the solution is sought. We denote the discrete Laplacian with mesh size $h_l$ by $\Delta^{h_l}$, with $\hl = \nref \hlp$. At each level we discretize the solution on a collection of node-centered cubic patches $\omli$, $R_l = \nref R_{l+1}$, and the corresponding discretized grids $\omlih$; the combined level $l$ grid is given by $\Omega^{l,h_l} \equiv \bigcup_\ibold \omlih$.  We also define, for each $\ibold$, localization regions 
$\omlia,\omliab$, and their discretizations  $\omliah,\omliabh$ 
$ 1 < \alpha < \beta$.  At level $0$ there is only one patch $\Omega^{0,h_0} $ at which the coarse solve of the method is performed, just as in the two-level algorithm.  We also impose a proper nesting condition: for $l=1...l_{max}$,
\begin{equation}
\label{eqn:propernesting2fdiscr}
\grow(\coarsen(\omliabh),b)\subset \Omega^{l-1,h_{l-1}}.
\end{equation}
The multilevel MLC comprises the following steps.

\noindent {\bf 1. Downward Pass - Initial Local Convolutions.}\\
Local convolutions are computed at levels $l=l_{max}, \dots , 1$.
\begin{eqnarray}
\label{eqn:mlcequationsfdiscr}
&&\phi^{\ibold,h_l}=G^{h_l}*\tilde{f}^{\ibold,h_l} \text{ on } \mathcal{G}(\omliah,N_{ref}b),
\end{eqnarray}
where the local right hand sides are defined as 
\begin{align*}
\tilde{f}^{\ibold,h_{l}}=&\sum \limits_{\ibold'} \Delta^{h_l}( \phi^{\ibold',h_{l+1} } ) |_{\mathcal{C}(\omliahp )} \\
+& \sum_{\ibold'}\Delta^{h_l}(G^{h_{l+1}}*\mathbb{P}(f^{\ibold',h_{l+1}}))|_{\mathcal{C}(\omliabhp \setminus \omliahp )}+ \tilde{\chi}_{\omlih} f^{h_l} 
\end{align*}
\begin{align*}
\tilde{\chi}_{\omlih}[\gbold] = \chi_{\omlih}[\gbold] - \sum \limits_{\substack{\ibold' = N_{ref} \ibold + \sbold  \\ 0 \leq \sbold_d \leq N_{ref}}} \chi_\omlihp[N_{ref}\gbold]
\end{align*}
\noindent {\bf 2. Global Coarse Solve }.
\begin{gather*}
\phi^{h_0}=G^{h_0}*\tilde{f}^{h_0} \text{ on }\Omega^{0,h_0}.
\end{gather*}

\noindent {\bf 3. Upward Pass - Local Interactions / Local Corrections for  $\boldsymbol{1} {\bf , \dots ,} \boldsymbol{l}_{\boldsymbol{max}}$}.\\
Starting from level 1, the following local Dirichlet problems are solved at levels $l=1,...,l_{max}$:
\begin{gather}
\label{eqn:last_step_multilevel}
\Delta^{h_l} \tilde{\phi}^{MLC,\ibold,h_l}=\tilde{f}^{\ibold,h_l} \text{ on } \omlih - \partial \omlih  ,\\
\tilde{\phi}^{MLC,\ibold,h_l}=\phi^{B,\ibold,h_l}\text{ on }\partial \omlih ,\nonumber\\
\tilde{\phi}^{MLC,l} = \tilde{\phi}^{MLC,\ibold,h_l} \text{ on } \omlih .\nonumber
\end{gather}
The Dirichlet boundary conditions are given by
\begin{equation}
\label{eqn:mlcmultilevel}
\phi^{B,\ibold,h_l}[\boldsymbol{g}]=\phi^{loc,l,\gbold}[\boldsymbol{g}]+\interp^{h_{l-1}}\left(\tilde{\phi}^{MLC,l-1} - \phi^{loc,l,\gbold} \right)(\boldsymbol{g}h_l) 
\end{equation}
Here the local potentials $\phi^{loc, \gbold, l}$ are given by:
\begin{equation}
\label{eqn:philocmlcmultilevel}
\phi^{loc,l,\gbold}[\gbold']=\sum \limits_{\ibold':\gbold  \in \Omega_{R_l,\ibold',\alpha}^{h^l}}\phi^{\ibold',h_l}[\gbold']
+\sum \limits_{\ibold':\boldsymbol{g} \in \Omega_{R_l,\ibold',\beta}^{h_l}\setminus \Omega_{R_l,\ibold',\alpha}^{h_l}}
(G^{h_l}*\mathbb{P}(f^{\ibold'}))[\gbold'].
\end{equation}
Finally, the Mehrstellen correction at all levels is applied as follows:

\begin{equation}
\label{eqn:mcapplicationfullmultipole}
\phi^{MLC,l}[\gbold]=\tilde{\phi}^{MLC,l}[\gbold]+ C_2 h_{l}^2f^{\ibold,{h_l}}[\gbold] \cma \boldsymbol{g} \in \Omega^{l,h_l}
\end{equation}

We do not have a complete error analysis for the above algorithm corresponding to that given in the two-level case. However, we can look at error analysis of the two-level algorithm, and determine the change in the error introduced there by replacing the coarse-grid convolution with $G^H$ with an MLC calculation.
We denote by:
\begin{itemize}
\item
$G^{MLC,S}(r)$ the two two-level semi-discrete method of local corrections approximation to $G*r$, with patch radius $S$;
\item
$N_1^S(r)(\xbold) \equiv \sum \limits_{\ibold: \xbold \notin \Omega_{S,\ibold, \beta}} h^q L^{q+2} (G*r^\ibold))(\xbold)$;
\item 
$N_2^S(r)(\xbold) \equiv \sum \limits_{\ibold: \xbold \in \Omega_{S,\ibold, \beta} \setminus  \Omega_{S,\ibold, \alpha}} h^q L^{q+2} (G*((\mathbb{I} - \mathbb{P})r^\ibold) )(\xbold) $; and
\item 
$N^S(r) = N_1^S(r) + N_2^S(r)$.
\end{itemize}
By \eqref{eqn:localcorrectionserror}, \eqref{eqn:SDC},
$G^H*(N^{R}(f))^H = (G*f)^H - \phi^H$
is the only quantity in the error in which convolution with $G^H$ appears. Given that, it is straightforward to assess the impact of replacing the convolution with $G^H$ in this expression with applying the MLC algorithm for a patch size $N_{ref} R$. To estimate this effect, we use a modified equation approach, in which the difference is approximated by $G*(N^R(f)) - G^{MLC,N_{ref} R}(N^R(f))$.
Applying the error estimate \eqref{eqn:SDC}, we obtain
\begin{align*}
G*(N^R(f)) - G^{MLC,N_{ref} R}(N^R(f)) = & N^{N_{ref}R}(N^R(f)) \\
= & N_1^{N_{ref}R}(N_1^R(f)) + N_1^{N_{ref}R}(N_2^R(f)) \\
+ & N_2^{N_{ref}R}(N_1^R(f)) + N_2^{N_{ref}R}(N_2^R(f))
\end{align*}
For this substitution to have an appropriately small impact, it is sufficient for the error to be comparable to or less than the error in the two-level algorithm. The sum of the first three terms meet this criterion -- the sum of first two terms is bounded by the max norm of the two-level error multiplied by $O(\beta^{-q})$, and the third term is bounded by $O(\alpha^{-q})$ times the max norm of the barrier error of the two-level algorithm. The final term, however, is problematic. In particular, the impact on the error of multiple applications of $\mathbb{I} - \mathbb{P}$ at increasing mesh spacings is far from clear. We will see evidence of this in the numerical results in Section \ref{sec:resOsc}, and will suggest a remedy that allows the error to be controlled.

\section{Computational Issues}

The analysis and demonstration of the performance of this algorithm will be the subject of a separate paper \cite{mlc2}, so we will just make a few high-level observations to justify the pursuit of this line of research. 
The largest contribution to the floating--point operation count in this method comes from the initial local discrete convolutions \eqref{eqn:initiallocalproblemsmlcII}. To compute these convolutions, we use a generalization of Hockney's domain-doubling algorithm \cite{hockneyMCP}, which we describe in the Appendix.  The floating point work per unknown for this step is $O(\alpha^3 log(n)), \alpha > 1$, where $n^3$ is the number of points per patch. The next-largest computation is that of the final Dirichlet solutions \eqref{eqn:last_step_fulldiscr}, performed using sine transforms, which is $O(log(n))$ per unknown. The floating point work associated with computing the Legendre expansions is small, with the convolutions of Legendre polynomials with the discrete Green's functions precomputed and stored. The memory overhead for storing these quantities scales like $O(\beta^3 n^3)$. However, there is one copy of these per processor, shared across multiple patches / multiple cores. Furthermore, they are only stored either on a sampled grid coarsened by $N_{ref}$, or on planar subsets corresponding to boundaries of patches, which reduces the memory overhead further.

The parallel implementation of this algorithm is via domain decomposition, with patches distributed to processors. For the choices of $\alpha$ and $\beta$ used in the results described below, this corresponds to a floating-point operation count about three times that of a corresponding multigrid algorithm for comparable accuracy. Roughly speaking, the communications costs, in terms of number of messages and overall volume of data moved, corresponds to that of a single multigrid V-cycle, plus the negligible costs of communicating a small number of Legendre expansion coefficients (20 per patch for the case $P=4$). This is to be compared to the eight multigrid V-cycles required to obtain a comparable level of accuracy. Current trends in the design of HPC processors based on low-power processor technologies indicate a rapid growth in the number of cores capable of performing floating-point operations on a processor, while the communications bandwith between processors, or between the processor and main memory, is growing much more slowly. In addition, most of the floating-point work is performed using FFTs on small patches on a single node, for which there are multiple opportunities for performance optimization. Thus the present algorithm is well-positioned to take advantage of these trends.  

\section{\large{Numerical Test Cases}}
\label{sec:examples}
We present in this section several examples that demonstrate the convergence properties of the MLC method described above. In all cases, we use as a measure of the solution error the max norm error of the potential, divided by max norm of the potential 
\begin{equation}
\frac{\|\phi^{MLC,h}-\phi\|_{\infty}}{||\phi\|_{\infty}}. \label{eqn:maxerr}
\end{equation}
For all cases, we set $n = 32$, so that $H/R = 1/4$.
We refer to the special case $\beta=\alpha$ (i.e. if the long-range potentials induced by the truncated Legendre expansions of local charges are ignored) as the \mlcOne{} method and to the general case $\alpha < \beta$ as the \mlcTwo{} method. It is not difficult to see that for \mlcOne{}, the estimate (\ref{eqn:finemesherrorMLC}) reduces to
\begin{equation}
\label{eqn:finemesherrorMLCaeqabar}
\phi^{MLC,h} - \phi=O(h^4)+O\left( h^2 \left(\frac{H}{R} \right)^{q_I-2}\frac{\|f\|_{\infty}}{\beta^{q_I-2}} \right)+O\left(\left(\frac{H}{R}\right)^q \frac{\|f\|_{\infty}}{\beta^q}\right)
\end{equation}
Increasing $\beta$ to reduce the barrier error in (\ref{eqn:finemesherrorMLCaeqabar}) substantially increases the per patch computational cost of the discrete convolution in the downward pass of the method. This is, in fact, the reason we replaced the local long-range potential values with the convolutions of the local Legendre expansions in Section \ref{sec:semiDiscreteAlgorithm}. 
\subsection{A Smooth Charge Distribution}
The first test case we are considering involves computing the potential induced by a smooth charge. The computational domain is the unit cube $\Omega=[0,1]^3$. The charge density is given by:
\begin{equation*}
f(\boldsymbol{x})=\left\{\begin{array}{r}(r-r^2)^4,r<1\\\\0,r\geq 1 \end{array} \right. , r=\frac{1}{R_o}\|\boldsymbol{x}-\boldsymbol{x}_o\|
\end{equation*}
and the support of the charge is a sphere of radius $R_o=\frac{1}{4}$, centered at point $\boldsymbol{x}_o=\frac{1}{2}\boldsymbol{1}$. The exact solution for this problem is given by:
\begin{equation*}
\phi(\boldsymbol{x})=R_o^2\left\{\begin{array}{r} \frac{r^6}{42}-\frac{r^7}{14}+\frac{r^8}{12}-\frac{2r^9}{45}+\frac{r^{10}}{110}-\frac{1}{1260},r<1 \\\\ -\frac{1}{2310r},r\geq 1 \end{array} \right.
\end{equation*}
and reduces to a pure monopole field for $r\geq 1$. 
\subsubsection{Two-Level Results}
\indent In Table \ref{tbl:mlcI} we present the fine mesh errors for the \mlcOne{} algorithm, with two levels, for mesh sizes $h=\frac{1}{256},\frac{1}{512},\frac{1}{1024}$ using the $L_{19}^h$ Mehrstellen Laplacian ($q=4$). We set $b = 2 \rightarrow q_I=6$ so that dependence of the interpolation error as a function of $\alpha,\beta$ matches that of the other error terms. 
For this problem, the errors in all three cases are so small that they are the barrier errors; each time we double $\beta$, the error goes down by roughly a factor of 16, as predicted by (\ref{eqn:finemesherrorMLCaeqabar}). 
\begin{table}[H]
\centerline{
\begin{tabular}{|l|l|l|l|l|l|}
\hline
 N          &  $\beta = 1.5$  & $\beta = 3.0$ & $\beta = 6.0$ \\
\hline 256  &  1.43756e-5      & 6.07186e-7     & 5.80288e-8     \\
\hline 512  &  1.29572e-5      & 4.32691e-7     & 2.67372e-8     \\
\hline 1024 &  1.27114e-5      & 4.01180e-7     & 2.44521e-8     \\
\hline
\end{tabular}
}
\caption{{\it 2-Level \mlcOne{}}: Scaled fine mesh maximum errors \eqref{eqn:maxerr} using the $L_{19}^h$ Mehrstellen Laplacian.}
\label{tbl:mlcI}
\end{table}
In Tables \ref{tbl:mlcIIaBar3L19} and \ref{tbl:mlcIIaBar6L19} we present fine mesh errors for the \mlcTwo{} algorithm, with $\alpha=1.5$, for $\beta=3$ and $\beta=6$, respectively, when refining both $h$ and $P$. As $h \rightarrow 0$, the error in this case approaches a barrier error for both the $P=1$ and $P=4$ cases at a rate of $O(h^2) - O(h^4)$, and those barrier errors correspond to the errors for the MLC-0 calculations with same corresponding values of $\beta$. For comparison, we also include the values of the error for the MLC-0 calculations with comparable computational costs, i.e. for $\beta = 1.5$. It is clear that for the negligible cost of adding the Legendre expansion, we obtain a decrease in the error by one-three orders of magnitude.
\begin{table}[H]
\centerline{
\begin{tabular}{|l|c|l|l|}
\hline
 N          &  \mlcOne{} ($\beta = 1.5$) & P=1           & P=4          \\
\hline 256  &  1.43756e-5                      & 4.35976e-6    & 1.63706e-6   \\
\hline 512  &  1.29572e-5                      & 1.43414e-6    & 4.58615e-7   \\ 
\hline 1024 &  1.27114e-5                      & 5.77475e-7    & 3.65246e-7   \\ 
\hline
\end{tabular}
}
\caption{{\it 2-Level \mlcTwo{}}: Scaled fine mesh maximum errors \eqref{eqn:maxerr} using $L_{19}^h$. For sufficiently small $h$ and $P=4$ nearly the same errors with the second column of Table \ref{tbl:mlcI} are obtained.}
\label{tbl:mlcIIaBar3L19}
\end{table}
\begin{table}[H]
\centerline{
\begin{tabular}{|l|c|l|l|l|}
\hline
 N          &  \mlcOne{} ($\beta = 1.5$) & P=1           & P=4        & P=5         \\
\hline 256  &  1.43756e-5                      & 4.05752e-6    & 1.45072e-6 & 1.68422e-6  \\
\hline 512  &  1.29572e-5                      & 1.12630e-6    & 1.04191e-7 & 4.49529e-8  \\
\hline 1024 &  1.27114e-5                      & 2.37651e-7    & 2.55964e-8 & 2.44951e-8  \\
\hline
\end{tabular}
}
\caption{{\it 2-Level \mlcTwo{}}: Scaled fine mesh maximum errors \eqref{eqn:maxerr} using $L_{19}^h$. Here $\alpha = 1.5,\beta=6$. For sufficiently small $h$ and high values of $P$ nearly the same errors with the third column of Table \ref{tbl:mlcI} are obtained.}
\label{tbl:mlcIIaBar6L19}
\end{table}

Next, we present the errors obtained by performing similar runs using the $L_{27}^h$ Mehrstellen Laplacian, for which $q=6$.  We set $b = 3 \rightarrow q_I=8$ so that dependence of the interpolation error as a function of $\alpha,\beta$ matches that of the other error terms.  In this case, the barrier error is $O(\beta^{-6})$; hence we expect that smaller values of the $\beta$ correction radius are required to obtain errors similar with those obtained with the $L_{19}^h$ difference operator. Since ${3^4}\approx{2^6},{6^4}\approx{3.25^6}$ we set $\beta=2,3.25$. 
First, in order to estimate the barrier values, we present the fine mesh errors for the \mlcOne{} method in  Table \ref{tbl:mlcIL27}  with $\beta=2,3.25$ using the $L_{27}^h$ operator. With those values of $\beta$; we expect comparable or smaller errors than those of the \mlcOne{} method with $\beta=3,6$ using the $L_{19}^h$ operator. This is the case, as is evident from a comparison with the error values of Table \ref{tbl:mlcI}. Furthermore, the barrier error as a function of $\beta$ decreases by more than the factor of $18.4 = (3.25/2)^6$ predicted by the analysis.

\begin{table}[H]
\centerline{
\begin{tabular}{|l|l|l|l|l|l|}
\hline
 N          &  $\beta = 2.0$      & $\beta = 3.25$ \\
\hline 256  &  1.25208e-7   & 4.11121e-8\\
\hline 512  &  1.14831e-7    & 4.92150e-9 \\
\hline 1024 &  1.01073e-7    & 3.11897e-9 \\
\hline
\end{tabular}
}
\caption{{\it 2-Level \mlcOne{}}: Scaled fine mesh maximum error \eqref{eqn:maxerr} using the $L_{27}^h$ Mehrstellen Laplacian. Compare with the second and third columns of Table \ref{tbl:mlcI}.}
\label{tbl:mlcIL27}
\end{table}
In Tables \ref{tbl:mlcIIaBar2L27} and \ref{tbl:mlcIIaBar3.3L27}, we present the errors for the MLC algorithm, for the cases $\beta=2,3.25$; $\alpha = 1.5$ for both cases. The $\beta = 2$ calculations reach the same barrier  errors as $h$ decreases. That is not the case for the $\beta = 3.25$ results in Table \ref{tbl:mlcIL27}, but that is not surprising -- the reduction of the barrier error by nearly an order of magnitude provides more headroom for $h$--convergence. However, we see that in Table \ref{tbl:mlcIIaBar3.3a1.75L27}, a slight increase of the inner correction radius to $\alpha=1.75$ allows us to reach the barrier error more rapidly. This is consistent with the error analysis, in that increasing $\alpha$ reduces the coefficient in front of the $O(h^P)$ error from truncating the Legendre expansion, from which we infer that the error from that source, rather than the error from the inner local convolution, is the dominant $h$-dependent error for this smooth example. 

\begin{table}[H]
\centerline{
\begin{tabular}{|l|l|l|l|l|l|l|}
\hline
 N          &  P=1           & P=4         \\
\hline 256  &  1.45293e-6    & 1.40270e-6  \\
\hline 512  &  5.20885e-7    & 1.89409e-7  \\               
\hline 1024 &  1.77613e-7    & 1.02341e-7  \\
\hline
\end{tabular}
}
\caption{{\it 2-Level \mlcTwo{}}: Scaled fine mesh maximum errors \eqref{eqn:maxerr} using $L_{27}^h$. Here $\alpha = 1.5,\beta=2 $. The $h \rightarrow 0$ errors are the same as the barrier errors in the first column of Table \ref{tbl:mlcIL27} }
\label{tbl:mlcIIaBar2L27}
\end{table}

\begin{table}[H]
\centerline{
\begin{tabular}{|l|l|l|l|l|l|l|}
\hline
 N          &   P=1           & P=4        & P=5         \\
\hline 256  &   1.40367e-6    & 1.47261e-6 & 1.63939e-6  \\
\hline 512  &   4.32214e-7    & 8.68274e-8 & 5.91126e-8  \\
\hline 1024 &   9.11841e-8    & 1.18905e-8 & 1.17441e-8  \\
\hline
\end{tabular}
}
\caption{{\it 2-Level \mlcTwo{}}: Scaled fine mesh maximum errors \eqref{eqn:maxerr} using $L_{27}^h$. Here $\alpha = 1.5,\beta=3.25$. The barrier errors are comparable with those using $L_{19}^h$ with $\beta=6$, (Table \ref{tbl:mlcIIaBar6L19}).}
\label{tbl:mlcIIaBar3.3L27}
\end{table}


\begin{table}[H]
\centerline{
\begin{tabular}{|r|l|l|l|l|l|l|l|}
\hline
 N          &   P=1           & P=4        & P=6        & P=9        \\
\hline 256  &   1.91470e-7    & 1.96490e-7 & 6.39837e-8 & 4.90745e-8 \\
\hline 512  &   5.42412e-8    & 9.16574e-9 & 5.99534e-9 & 6.02843e-9 \\   
\hline 1024 &   1.40428e-8    & 2.79547e-9 &            &            \\ 
\hline
\end{tabular}
}
\caption{{\it 2-Level \mlcTwo{}}: Scaled fine mesh maximum errors \eqref{eqn:maxerr} using $L_{27}^h$. Here $\alpha = 1.75,\beta=3.25$. Compare with the second column of Table \ref{tbl:mlcIL27}. A high polynomial degree is required to attain it for $h=\frac{1}{256}$.}
\label{tbl:mlcIIaBar3.3a1.75L27}
\end{table}

\subsubsection{Three-Level Results}
\label{sec:3Lsmooth}
We next present similar results using the multilevel MLC algorithm of Section \ref{sec:multilevel} with three levels. Since we have demonstrated a clear advantage to using the 27-point stencil, in the remaining studies we will restrict our attention to that operator.
In Table \ref{tbl:mlcIL27ThreeLevels} we show the barrier fine mesh errors obtained using the \mlcOne{} method for $\beta=2,3.25$. The errors for $\beta=3.25$ are more than 18.4 times smaller than the errors for $\beta=2$ and are nearly the same to the 2-level method errors (Table \ref{tbl:mlcIL27}). As predicted by the error analysis in Section \ref{sec:multilevel} the error of MLC-0 is insensitive to the number of levels. 

\begin{table}[H]
\centerline{
\begin{tabular}{|l|l|l|l|l|l|}
\hline
 N          &  $\beta = 2.0$       & $\beta = 3.25$     \\
\hline 512  &  1.30594e-7   & 4.86092e-9 \\
\hline 1024 &  1.90632e-7  & 3.92874e-9 \\
\hline
\end{tabular}
}
\caption{{\it 3-Level \mlcOne{}}: Scaled fine mesh maximum errors \eqref{eqn:maxerr} using the $L_{27}^h$ Mehrstellen Laplacian. Compare with Table \ref{tbl:mlcIL27} that contains the 2-level results.}
\label{tbl:mlcIL27ThreeLevels}
\end{table}

In Table \ref{tbl:mlcII3levelaBar3.25L27.smooth} the errors obtained with the 3-level \mlcTwo{} method are shown using $\alpha = 1.75$,  $\beta=3.25$. Unlike the two-level results, the $P=4$ errors are significantly poorer than the MLC-0 errors.  For example, we recover the barrier errors only for $N=4096$, as opposed to the $N=512$ results for MLC-0. We can improve matters somewhat by increasing $P$, but even for this very smooth problem, we do not get close to the barrier errors until $N=2048$. This is consistent with the analysis in Section \ref{sec:multilevel}, and indicates that using higher values of $P$ does not solve the problem. We will propose a different solution in Section \ref{sec:resOsc}.

\begin{table}[H]
\centerline{
\begin{tabular}{|r|l|l|l|l|}
\hline
        N, level     & P=1       & P=4       & P=6        & P=8       \\
\hline  512  l=0     & 1.4509e-7 & 1.1379e-7 & 5.8886e-8  & 4.2602e-8 \\
\hline       l=1     & 4.9396e-7 & 1.0594e-6 & 1.0990e-6  & 3.0059e-7 \\
\hline       l=2     & 5.2600e-7 & 1.0782e-6 & 1.1101e-6  & 1.4926e-7 \\
\hline 1024  l=0     & 1.1143e-7 & 1.9032e-8 & 9.5197e-9  & 4.1018e-9 \\
\hline       l=1     & 2.2539e-7 & 1.6461e-7 & 9.9491e-8  & 2.3381e-8 \\
\hline       l=2     & 2.3665e-7 & 1.6596e-7 & 9.9989e-8  & 2.3381e-8 \\
\hline 2048  l=0     & 3.8485e-8 & 5.7487e-9 & 5.0311e-9  & \\
\hline       l=1     & 5.9923e-8 & 1.0143e-8 & 5.9864e-9  & \\
\hline       l=2     & 6.1989e-8 & 1.0276e-8 & 6.0168e-9  & \\
\hline 4096  l=0     & 1.3028e-8 & 5.1364e-9 &            & \\
\hline       l=1     & 1.6487e-8 & 5.2147e-9 &            & \\
\hline       l=2     & 1.6861e-8 & 5.2621e-9 &            & \\
\hline               
\end{tabular}
}
\caption{{\it 3-Level \mlcTwo{}}: Scaled maximum errors \eqref{eqn:maxerr} at all levels using $L_{27}^h$. Here $\alpha=1.75,\beta=3.25$. Compare with the second column of Table \ref{tbl:mlcIL27ThreeLevels}.}
\label{tbl:mlcII3levelaBar3.25L27.smooth}
\end{table}

\subsection{An Oscillatory Charge Test Case} \label{sec:resOsc}
We further consider a case of three oscillatory charges that has been previously studied in \cite{mlc}. The computational domain is again the unit cube $\Omega=[0,1]^3$. Here we define a local charge density, whose support is a sphere of radius $R_o$ centered at point $\boldsymbol{x}_o$, by:
\begin{equation}
\label{eqn:oscillatorycharge}
f_{\boldsymbol{x}_o}(\boldsymbol{x})=\left\{\begin{array}{r}\frac{1}{R_o^3}(r-r^2)^2\sin^2(\frac{\gamma}{2}r),r<1\\\\0,r\geq 1 \end{array} \right. , r=\frac{1}{R_o}\|\boldsymbol{x}-\boldsymbol{x}_o\|,\gamma=4\mu\pi,\mu=7
\end{equation}
The exact solution associated with this charge density is given by:
\begin{equation*}
\phi_{\boldsymbol{x}_o}(\boldsymbol{x})=\frac{1}{R_o}\left\{\begin{array}{l} -\frac{1}{120} - \frac{6}{\gamma^4}\:\:\:\:\:\:\:\:\:\:\:\:\:\:\:\:\:\:\:\:\:\:\:\:\:\:\:\:\:\:\:\:\:\:\:\:\:\:\:\:\:\:\:\:\:\:\:\:\:\:\:\:\:\:\:\:\:\:\:\:\:\:\:\:\:\:\:\:\:\:\:\:\:\:\:\:\:\:\:\:\:\:\:\:\:\:\:\:\:\:\:\:\:\:\:\:\:\:\:\:,\:\:r=0 \\\\
\frac{r^6}{84} - \frac{r^5}{30} + \frac{r^4}{40} + \frac{60}{\gamma^6} - \frac{9}{\gamma^4} - \frac{1}{120} + \frac{120}{\gamma^6r} \\\\
 + \left(-\frac{120}{\gamma^6r} - \frac{9}{\gamma^4} + \frac{300}{\gamma^6} + \frac{36r}{\gamma^4} + \frac{r^2}{2\gamma^2} - \frac{30r^2}{\gamma^4} - \frac{r^3}{\gamma^2} + \frac{r^4}{2\gamma^2}\right) \cos(\gamma r)  \\\\
 + \left(\frac{12}{\gamma^5r} - \frac{360}{\gamma^7r} - \frac{96}{\gamma^5} + \frac{120r}{\gamma^5} - \frac{3r}{\gamma^3} + \frac{8r^2}{\gamma^3} - \frac{5r^3}{\gamma^3}\right)  sin(\gamma r)\:\:\:\:\:\:\:\:\:\:\:\:\:\:\:\:\:\:,\:\:r<1 \\\\ 
\left(-\frac{1}{210} - \frac{12}{\gamma^4} + \frac{360}{\gamma^6}\right)\frac{1}{r}\:\:\:\:\:\:\:\:\:\:\:\:\:\:\:\:\:\:\:\:\:\:\:\:\:\:\:\:\:\:\:\:\:\:\:\:\:\:\:\:\:\:\:\:\:\:\:\:\:\:\:\:\:\:\:\:\:\:\:\:\:\:\:\:\:\:\:\:\:\:\:\:\:\:\:\:\:\:\:\:,\:\:r\geq 1 \end{array} \right.
\end{equation*}
and is a pure monopole for $\:r\geq 1\:$. For our test case we consider three charges of the form (\ref{eqn:oscillatorycharge}), of radius $R_o=\frac{5}{100}$, centered at points $\boldsymbol{c}_1=\left(\frac{3}{16},\frac{7}{16},\frac{13}{16} \right)$, $\boldsymbol{c}_2=\left(\frac{7}{16},\frac{13}{16},\frac{3}{16} \right)$ and $\boldsymbol{c}_3=\left(\frac{13}{16},\frac{3}{16},\frac{7}{16} \right)$. The total charge and total potential are given via linear superposition by:
\begin{eqnarray*}
&&f(\boldsymbol{x})=f_{\boldsymbol{c}_1}(\boldsymbol{x})+f_{\boldsymbol{c}_2}(\boldsymbol{x})+f_{\boldsymbol{c}_3}(\boldsymbol{x})\\
&&\phi(\boldsymbol{x})=\phi_{\boldsymbol{c}_1}(\boldsymbol{x})+\phi_{\boldsymbol{c}_2}(\boldsymbol{x})+\phi_{\boldsymbol{c}_3}(\boldsymbol{x})
\end{eqnarray*}

We first present the results using three levels (Table \ref{tbl:mlcI3levela3.25L27b3}) and four levels (Table \ref{tbl:mlcI4levela3.25L27b2}) using \mlcOne{}. The primary features of the convergence properties of the solution are that the errors are nearly uniform as a function of level, and are the same in both the three and four level cases. There is some indication of slowing down of the convergence rate on the finest two levels, but the convergence is still faster than $O(h^2)$
\begin{table}[H]
\centerline{
\begin{tabular}{|r|l|}
\hline
                 N, level    & error     \\         
\hline           2048 l=0    & 9.59918e-7 \\
\hline                l=1    & 1.00600e-6 \\
\hline                l=2    & 1.04402e-6 \\
\hline           4096 l=0    & 5.82005e-8 \\
\hline                l=1    & 6.47409e-8 \\
\hline                l=2    & 6.71067e-8 \\
\hline           8192 l=0    & 8.42867e-9 \\ 
\hline                l=1    & 8.42867e-9 \\
\hline                l=2    & 8.44657e-9 \\
\hline               
\end{tabular}
}
\caption{{\it 3-Level \mlcOne{}}: Scaled maximum errors \eqref{eqn:maxerr}, $\beta=3.25$.}
\label{tbl:mlcI3levela3.25L27b3}
\end{table}

\begin{table}[H]
\centerline{
\begin{tabular}{|r|l|}
\hline
                 N, level    &  error             \\         
\hline           2048 l=0      & 1.03645e-7    \\                           
\hline                l=1    & 9.59723e-7    \\                           
\hline                l=2    & 1.00621e-6   \\                           
\hline                l=3    & 1.04423e-6   \\                           
\hline           4096 l=0   & 2.93837e-8    \\                          
\hline                l=1   & 5.84863e-8   \\                          
\hline                l=2    & 6.50247e-8  \\                           
\hline                l=3    & 6.73912e-8    \\                           
\hline           8192 l=0      & 7.84890e-9    \\                           
\hline                l=1     & 8.78853e-9   \\                        
\hline                l=2      & 8.78853e-9    \\                          
\hline                l=3    &  8.79911e-9   \\                           
\hline               
\end{tabular}
}
\caption{{\it 4-Level \mlcOne{}}: Scaled maximum errors \eqref{eqn:maxerr}, $\beta=3.25$.}
\label{tbl:mlcI4levela3.25L27b2}
\end{table}

In the \mlcTwo{} convergence results in Table \ref{tbl:mlcII4levelaBar3.25L27b3}, we see substantial deviations from the \mlcOne{} convergence results. The error shows no consistent behavior as a function of resolution, and in fact is worse at the finest resolution (N = 8192) in Table \ref{tbl:mlcII4levelaBar3.25L27b3} than it is at the N=2048 resolution in Table \ref{tbl:mlcI4levela3.25L27b2}. We see no analogous problems in the \mlcOne{} calculations. Examining the error analysis in Section \ref{sec:multilevel}, we identified the terms in a three-level calculation that might lead to problems.   Even in the smooth example above, it is clear that the increasing $P$ does not have sufficient impact to solve this problem.  A different approach, suggested by the form of the error, is to reduce the difference $\beta - \alpha$ at coarser levels. In fact, there is likely a mechanism for defining a systematic strategy for doing this, since $(\mathbb{I} - \mathbb{P}) f^\ibold$ is easily computed. We defer that to later work. For the moment, we demonstrate this by setting $\alpha$ as an empiricially-determined slowly decreasing function of level, holding $\beta$ fixed (Table \ref{tbl:mlcII4levelaBar3.25L27b3VaryingAlpha}). We see that we can recover close to the errors in the \mlcOne{} calculation. In addition, the cost of increasing $\alpha$ slightly at coarser levels has a small impact of the overall cost of a multiresolution calculation, since these are applied to calculations at the coarser resolutions, which remain a small fraction of the overall cost of the method, even with the increased values of $\alpha$.

\begin{table}[H]
\centerline{
\begin{tabular}{|r|l|l|}
\hline
                 N, level                      & P=1        & P=4        \\
\hline           2048 l=0                      & 1.49110e-6 & 1.22696e-6 \\
\hline                l=1                      & 2.68289e-6 & 2.19207e-6 \\
\hline                l=2                      & 2.24874e-6 & 2.85969e-6 \\
\hline                l=3                      & 2.74856e-6 & 3.40770e-6 \\
\hline           4096 l=0                      & 3.26018e-7 & 2.37540e-7 \\
\hline                l=1                      & 3.78691e-7 & 2.03681e-6 \\
\hline                l=2                      & 6.21459e-7 & 6.97742e-7 \\
\hline                l=3                      & 6.42162e-7 & 7.16144e-7 \\
\hline           8192 l=0                      & 3.26374e-7 & 3.30828e-7 \\
\hline                l=1                      & 2.04758e-6 & 2.39355e-6 \\
\hline                l=2                      & 2.13626e-6 & 2.41492e-6 \\
\hline                l=3                      & 2.14514e-6 & 2.41812e-6 \\
\hline               
\end{tabular}
}
\caption{{\it 4-Level \mlcTwo{}}: Scaled maximum errors \eqref{eqn:maxerr}. Here $\:\alpha=1.75,\:\beta=3.25$.}
\label{tbl:mlcII4levelaBar3.25L27b3}
\end{table}

\begin{table}[H]
\centerline{
\begin{tabular}{|r|l|l|l|}
\hline
                 N, level   &   $\alpha_l$ & P=1        & P=4        \\
\hline           2048 l=0   &      -       & 1.53798e-7 & 1.48467e-7 \\
\hline                l=1   &    2.25      & 9.08569e-7 & 9.32892e-7 \\ 
\hline                l=2   &    2         & 1.05388e-6 & 1.14004e-6 \\ 
\hline                l=3   &    1.75      & 1.22595e-6 & 1.35075e-6 \\ 
\hline           4096 l=0   &      -       & 3.43626e-8 & 2.83645e-8 \\
\hline                l=1   &    2.25      & 6.15251e-8 & 1.00481e-7 \\
\hline                l=2   &    2.25      & 6.44820e-8 & 7.59802e-8 \\
\hline                l=3   &    2         & 7.29285e-8 & 8.29560e-8 \\
\hline           8192 l=0   &      -       & 7.17091e-9  & 7.48667e-9 \\
\hline                l=1   &    2.75      & 1.40237e-8  & 1.87350e-8 \\
\hline                l=2   &    2.75      & 1.43363e-8  & 1.92898e-8 \\
\hline                l=3   &    2         & 9.03511e-9  & 1.00333e-8 \\
\hline               
\end{tabular}
}
\caption{{\it 4-Level \mlcTwo{}}: Scaled maximum errors \eqref{eqn:maxerr} using $\:L_{27}^h\:$ with higher values of $\:\alpha\:$ at intermediate levels. Here $\:\beta=3.25$. Compare with Table \ref{tbl:mlcI4levela3.25L27b2}.}
\label{tbl:mlcII4levelaBar3.25L27b3VaryingAlpha}
\end{table}

\section{Conclusions}
We have presented a domain decomposition method for the numerical solution of Poisson's equation with infinite domain boundary conditions in three dimensions on a nested hierarchy of structured grids. The method is an extension of Anderson's Method of Local Corrections for particles \cite{anderson} to gridded data and generalizes the scheme of McCorquodale, et al. \cite{mlc}. 
In the present method, local potentials are computed as volume potentials of local charges up to an inner localization radius, combined with volume potentials induced by order $P-1$ truncated Legendre expansions of the local charges up to an outer localization radius. The remaining global coupling is represented using a coarse-grid version of the same representation. This generalizes the method in \cite{mlc}, which corresponds to the $P=1$ special case in the current method. Also, in \cite{mlc} the local potentials were computed by means of the James-Lackner representation \cite{james,lackner} of infinite--domain boundary conditions. In the present work, this is replaced by a representation using discrete convolution operators, which can be computed efficiently using FFTs via Hockney's algorithm. This approach eliminates the complicated quadratures that are necessary for the extension of the James-Lackner algorithm to three dimensions, while the FFT-based approach leads to compact compute kernels that can be highly optimized. The resulting algorithm is well-suited for high performance on HPC computing platforms made up of multicore processors; in \cite{mlc2}, we will present a systematic study of the performance and scaling of the algorithm on such systems.

In this paper, we have focused primarily on the analytical foundations of the MLC method and have provided a detailed error analysis. The errors are of the form $O(h^P) + O(h^4) + O(h^2 \beta^{-q} ) + O(\beta^{-q})$, where $h$ is the mesh spacing, $\beta$ is the nondimensionalized outer localization radius which is independent of $h$, and $q$ is the order of accuracy of the Mehrstellen operator on harmonic functions. Numerical experiments indicate that the observed convergence behavior of the method is consistent with the analysis. For computationally practical values of the localization radius, and using the 27-point Mehrstellen operator (for which $q = 6$), the barrier error corresponds to relative solution error norms of $10^{-8} - 10^{-9}$. While the $\beta^{-q}$ term looks like an $O(1)$ error relative to the mesh spacing $h$, it is better to think of it as a separate discretization parameter that governs the accuracy of the representation of the nonlocal coupling. Doubling $\beta$ decreases the error by a factor of $2^{-q}$, analogous to the impact of halving $h$. 

For the two-level algorithm, the results indicate that, for a given choice of the Mehrstellen operator, the two localization radii, and for $P=4$, the method converges at a rates  in the range $O(h^4)$--$O(h^2)$, until the error reaches the barrier, i.e. consistent with the error analysis. We have also defined and implemented the extension to more than two levels, following the approach in \cite{mlc}. A preliminary analysis of that algorithm indicates the need to control errors at coarser levels coming from the field induced between the inner and outer localization radii by the truncation of the Legendre expansion. The analysis suggests that these might be controlled by increasing the inner localization radius $\alpha$ at coarser levels. The numerical examples indicate that the problem is real, and that the proposed solution represents a viable approach. More generally, an important question that needs to be addressed is turning the error analysis in this work into practical strategies for choosing discretization parameters. For example, what are the tradeoffs between decreasing $\beta - \alpha$ and decreasing $h$ in order to improve the accuracy of a calculation, versus the cost of doing each? We will address these issues in \cite{mlc2}.

There are various possible ways to extend the present work. Perhaps most straightforward are extensions to finite--volume discretizations and the implementation of other boundary conditions on rectangular domains (including periodic boundary conditions) using a method--of--images approach. Another possibility would be to apply even higher--order Mehrstellen discretizations of the Laplacian to see whether it results in smaller values of the barrier errors than those reported in this work. As was seen in Section \ref{sec:examples}, the $\:L_{27}^h\:$ ($q = 6$) Mehrstellen Laplacian leads to comparable barrier errors to those obtained using the $\:L_{19}^h\:$ ($q=4$) stencil, but using smaller localization radii, in a manner consistent with the $O(\beta^{-q})$ scaling of that error. It is possible to derive Mehrstellen stencils for which $q = 10$, with the stencil contained in a $5 \times 5 \times 5$ block around the evaluation point. This leads to only a modest increase in the computational cost and complexity: for example, the per-patch computational cost of the most compute-intensive component of the algorithm -- the local discrete convolutions -- does not depend on the size of the stencil. Finally, it would be interesting to investigate extensions of this method to other elliptic problems in mathematical physics employing different Green's functions and high-order discretizations of the associated differential operators. The error analysis of the method as extended to other kernels should be essentially the same with what is discussed in the present study. Moreover, Hockney's algorithm is kernel-independent and can be readily applied with minor modifications.
More generally, the present work uses some detailed analytic tools for understanding the discrete potential theory on locally--structured grids associated with the combination of finite difference localization in \cite{mayo} and the local interactions / local corrections construction underlying \cite{anderson}. It would be interesting to go back to the original MLC method for particles and to other particle-grid methods, such as particle-in-cell and immersed boundary methods, and apply these tools to better understand the error properties of these methods.

\appendix
\section{{Appendix}}
\subsection{$L_{19}^h\:$ and $\:L_{27}^h\:$ Mehrstellen Discretizations of the Laplacian}
\label{sec:L19L27}

The stencil coefficients for the $\:L_{19}^h\:$ and $\:L_{27}^h\:$ Mehrstellen Laplacians are $\:a_{\boldsymbol{j}}=\frac{1}{h^2}b_{|\boldsymbol{j}|}$, where $\:|\boldsymbol{j}|\:$ is the number of non-zero components of $\:\boldsymbol{j}$ and $\:b_{k}\:$ are defined as:
\begin{eqnarray*}
b_0=-4,&&\:\:b_1=\frac{1}{3},\:\:\:\:\:\:\:\:b_2=\frac{1}{6},\:\:\:\:\:\:\:\:\:\:\:b_3=0,\:\:\:\:\: \text{19-point stencil}\\
b_0=-\frac{64}{15},&&\:\:b_1=\frac{7}{15},\:\:\:\:\:\:b_2=\frac{1}{10},\:\:\:\:\:\:\:\:b_3=\frac{1}{30},\:\: \text{27-point stencil}\\
\end{eqnarray*}

The corresponding expressions  for the truncation errors $\tau_{19}^h$, $\tau_{27}^h$ for $\:L_{19}^h\:$, $\:L_{27}^h\:$, are given by:
\begin{eqnarray}
\label{eqn:truncation19}
\nonumber \tau_{19}^h(\phi)= \frac{h^2}{12}(\Delta(\Delta \phi))+
h^4 L^{(6)}(\phi)  + O(h^6)
\end{eqnarray}
and 
\begin{eqnarray}
\label{eqn:truncation27}
\nonumber \tau_{27}^h(\phi)=\frac{h^2}{12}(\Delta(\Delta \phi))+\frac{h^4}{360}\left(\left(\Delta^2+2\left(\frac{\partial^4}{\partial x^2\partial y^2}+\frac{\partial^4}{\partial y^2\partial z^2}+\frac{\partial^4}{\partial z^2\partial x^2} \right)\right)(\Delta \phi)\right)+h^6 L^{(8)}(\phi)  + O(h^8)
\end{eqnarray}
where the $L^{(q)}$'s are homogeneous constant--coefficient $q^{th}$-order differential operators.

We need to compute an approximation to the discrete Green's function \eqref{eqn:dgf} for the 19-point and 27-point operators, restricted to a domain of the form $\:D=[-n,\:n]^3$. We do this by solving the following inhomogeneous Dirichlet problem on a larger domain $D_\zeta = [-\zeta n, \zeta n]^3$ .
\begin{align*}
\nonumber (L^{h=1} G^{h=1})[\gbold] = &\delta_{\boldsymbol{0}}[\gbold] \:\text{for}\: \gbold \in \mathcal{G}(D_\zeta,-1), \\
G^{h=1}[\gbold] =  &G(\gbold) \text{ for } \gbold \in D_\zeta - \mathcal{G}(D_\zeta,-1).
\end{align*}
where $G = G(\xbold)$ is the Green's function \eqref{eqn:potential}, and $L^h$ is either the 19-point or 27-point operator. Then our approximation to $G^{h=1}$ on $D$ is the solution computed on $D_\zeta$, restricted to $D$. To compute this solution, we put the inhomogeneous boundary condition into residual-correction form, and solve the resulting homogeneous Dirichlet problem using the discrete sine transform. 
The error estimate \eqref{eqn:msCorr4} applied here implies that the error in replacing the correct discrete boundary conditions with those of the exact Green's function scales like $O({(\zeta n)}^{-4})$ in max norm. In the calculations presented here, we computed $G^{h=1}$ using $n \ge 128$ and $\zeta = 2$, leading to at least 10 digits of accuracy for $G^{h=1}$. 

\subsection{Hockney's Method for Fast Evaluation of Discrete Convolutions}

Hockney (\cite{hockneyMCP},p.180--181; see also \cite{eastwood}) observed that discrete convolutions with one of the functions having support on a bounded domain in $\mathbb{Z}^\Dim $, and evaluated on a bounded domain, can be computed exactly in terms of discrete Fourier transforms. For completeness, we describe that method. We show this first for the case $\Dim = 1$, and state the general result for any number of dimensions. Given $\Psi , f : \mathbb{Z} \rightarrow \mathbb{R}$, $supp(f) \subseteq [0,b]$, we want to compute 
\begin{gather}
(\Psi * f) [i] = (f * \Psi )[i] = \sum \limits_{j \in \mathbb{Z}} f[i - j]  \Psi[j], i \in [0,n]
\end{gather}
First, we observe that the infinite sum can be replaced by a finite sum.
\begin{gather} \label{eqn:hockney1}
\sum \limits_{j \in \mathbb{Z}} f[i - j]  \Psi[j] = \sum \limits_{j = -b'}^n f[i - j]  \Psi[j], i \in [0,n]
\end{gather}
for any $b' \geq b$. Second, we observe that $\Psi$, $f$ can be replaced in \eqref{eqn:hockney1} by periodic extensions of those functions restricted to the interval $[-b',n]$.
\begin{gather}
\sum \limits_{j = -b'}^n f[i - j]  \Psi[j] = \sum \limits_{j = -b'}^n \tilde{f}[i - j]  \tilde{\Psi}[j], i \in [0,n] \label{eqn:hockney2}\\
\tilde{f}[l] , \tilde{\Psi}[l] \equiv f[l_{mod}],\Psi[l_{mod}] \text{ , } l_{mod} = mod(l + b', (n + b' + 1)) - b'. \nonumber
\end{gather}
Finally, we express the periodic convolution in \eqref{eqn:hockney2} in terms of discrete Fourier transforms. 
\begin{gather}
\sum \limits_{j = -b'}^n \tilde{f}[i - j]  \tilde{\Psi}[j] = \mathcal{F}^{-1} (\mathcal{F}(\tilde{\Psi}) \cdot \mathcal{F} (\tilde{f}))[i],
\end{gather}
where $\mathcal{F}$, $\mathcal{F}^{-1}$ are the discrete complex Fourier transform and its inverse on the interval $[-b',n] \subset \mathbb{Z}$. 

This generalizes to rectangular domains in any number of dimensions. For example, for cubic domains,
given $\Psi , f : \mathbb{Z}^\Dim \rightarrow \mathbb{R}^\Dim$, $supp(f) \subseteq [0,b]^\Dim$, 
\begin{gather}
\sum \limits_{\jbold \in \mathbb{Z}^\Dim} \Psi[\ibold - \jbold] f[\jbold] = \mathcal{F}^{-1}(\mathcal{F}(\tilde{\Psi}) \cdot \mathcal{F}( \tilde{f}) )[\ibold] , \ibold \in [0,n]^\Dim  \\
\tilde{f}[\lbold] , \tilde{\Psi}[\lbold] \equiv f[\lbold_{mod}],\Psi[\lbold_{mod}], \\ (\lbold_{mod})_d = mod((\lbold)_d + b', (n + b' + 1)) - b' , d = 0, \dots \Dim-1,
\end{gather}
where $b' \geq b$ and $\mathcal{F}$, $\mathcal{F}^{-1}$ are the complex discrete Fourier transform and its inverse on the cube $[-b',n]^\Dim \subset \mathbb{Z}^\Dim$.
In practice, this is efficient for a broad range of $(b,n)$ since we can choose $b'$ so that the radices of the FFTs are highly composite, with the size of the problem changing by only a small amount. In the case where $b=n$, the length of the domain doubles in each direction, hence this is often referred to as Hockney's domain-doubling algorithm. However, in the present application, we want to use the more general case, since the size of the support of the localized charge distributions and the size of the grid on which the local fields are defined differ by a significant amount.
\section*{Acknowledgments}
The authors would like to thank Brian Van Straalen and Peter McCorquodale for a number of helpful discussions. 
This research was supported at the Lawrence Berkeley National Laboratory by the Office of Advanced Scientific Computing Research of the U.S. Department of Energy under Contract No. DE-AC02-05CH11231 and at the National Energy Research Scientific Computing Center by the DOE Petascale Initiative in Computational Science and Engineering.


\end{document}